\documentclass[10pt]{amsart}

\usepackage{amsmath}
\usepackage{mathtools}
\usepackage{enumerate}
\usepackage{amssymb}
\usepackage{tensor}
\usepackage{etoolbox}
\usepackage{dirtytalk}
\usepackage{comment}
\usepackage{enumitem}
\usepackage{hyperref}

\usepackage{tikz}
\usetikzlibrary{shapes.geometric}
\usetikzlibrary{cd}
\usetikzlibrary{calc}
\tikzcdset{every label/.append style = {font = \small}}

\newcommand{\A}{\mathcal{A}}
\newcommand{\B}{\mathcal{B}}
\newcommand{\C}{\mathcal{C}}
\newcommand{\D}{\mathcal{D}}

\newcommand{\F}{\mathcal{F}}

\newcommand{\I}{\mathcal{I}}
\newcommand{\J}{\mathcal{J}}
\newcommand{\K}{\mathcal{K}}
\newcommand{\M}{\mathcal{M}}

\newcommand{\Z}{\mathbb{Z}}

\newcommand{\Sc}{\mathcal{S}}

\newcommand{\Q}{\mathcal{Q}}

\newcommand{\cZ}{\mathcal{Z}}

\newcommand{\e}{\epsilon}

\newcommand{\simp}[2]{e_{#1#2}}
\newcommand{\Nim}{\text{\fontfamily{cmss}\selectfont N}}
\newcommand{\Enc}{\text{\fontfamily{cmss}\selectfont E}}
\DeclareMathOperator{\Gr}{K}

\newcommand{\opalg}{A_{\text{op}}}
\newcommand{\clalg}{A_{\text{cl}}}

\newcommand{\RC}{\mathcal{RC}}

\newcommand{\RTC}{\mathcal{RTC}}

\newcommand{\TC}{\mathcal{TC}}
\newcommand{\TM}{\mathcal{TM}}

\DeclareMathOperator{\Obj}{Obj}
\DeclareMathOperator{\Irr}{Irr}
\DeclareMathOperator{\End}{End}
\DeclareMathOperator{\Mod}{Mod}
\DeclareMathOperator{\Hom}{Hom}
\DeclareMathOperator{\Aut}{Aut}
\DeclareMathOperator{\my-hom}{hom \!}
\DeclareMathOperator{\Nat}{Nat}
\DeclareMathOperator{\Spec}{Spec}

\DeclareMathOperator{\Tr}{{\normalfont Tr}}
\DeclareMathOperator{\tr}{{\normalfont tr}}
\DeclareMathOperator{\id}{id}
\DeclareMathOperator{\Vect}{\underline{Vect}}
\DeclareMathOperator{\cre}{cre}
\DeclareMathOperator{\ann}{ann}
\DeclareMathOperator{\Ima}{Im}

\newcommand{\rtcs}[2]{\epsilon_{#1}^{#2}}
\newcommand{\CC}{\C \boxtimes \overline{\C}}
\newcommand{\Blank}{\, \vcenter{\hbox{\tiny$\bullet$}} \, }
\newcommand{\fld}{\mathbb{K}}

\newcommand{\tid}{\normalfont \textbf{1}}

\newcommand{\vprod}{\cdot}

\newcommand{\tl}[1][\beta]{\text{\normalfont TL}}
\newcommand{\Rtl}[1][\beta]{\mathcal{R}\text{\normalfont TL}}
\newcommand{\tlred}[1][\beta]{\text{\normalfont TL}^\text{red}}
\newcommand{\Rtlred}[1][\beta]{\mathcal{R}\text{\normalfont TL}^\text{red}}

\newcommand{\<}{\langle}
\renewcommand{\>}{\rangle}

\newcommand{\Sharp}{^\sharp}
\newcommand{\Flat}{^\flat}
\newcommand{\Natural}{^\natural}

\newcommand{\defeq}{\vcentcolon=}

\newcommand{\mult}{[d,\chi_{\M}]}

\theoremstyle{plain}
\newtheorem{THM}{Theorem}[section]
\newtheorem{PROP}[THM]{Proposition}
\newtheorem{LEMMA}[THM]{Lemma}

\newtheorem{COR}[THM]{Corollary}
\theoremstyle{definition}
\newtheorem{REM}[THM]{Remark}
\newtheorem{EX}[THM]{Example}
\newtheorem{DEF}[THM]{Definition}
\newtheorem*{DEF*}{Definition}

\def\mystrut(#1,#2){\vrule height #1pt depth #2pt width 0pt}
\title{Modular invariants and NIM-reps}
\author{Alastair King}
\address{
  \parbox{\linewidth}{
Department of Mathematical Sciences, University of Bath,\\ Bath, BA2 7AY, United Kingdom}
}

\author{Leonard Hardiman}
\address{
  \parbox{\linewidth}{
CMS, École Polytechnique Fédérale de Lausanne,\\
EPFL, Station 4, 1015 Lausanne, Suisse}
}

\email{leonard.hardiman@epfl.ch}

\begin{document}

\begin{abstract}
  Given a pivotal module category over a spherical fusion category, we introduce the encircling module, a module over the fusion algebra defined using the pivotal structure, and prove that it is isomorphic to the NIM-rep as a fusion algebra module. When applied to the $\TM$ realisation of the modular invariant partition function~\cite{Hardiman:2019utc}, this identifies the NIM-rep multiplicities with the entries $Z_{I,I^{\vee}}$ of the modular invariant (its diagonal relative to the charge-conjugation invariant, which is the one given by the regular module category), providing a categorical generalisation of Boeckenhauer, Evans and Kawahigashi's results~\cite{Boeckenhauer1999ChiralSO}. We also show that for indecomposable module categories the dimension condition on $\TM$ required for modular invariance is automatically satisfied, and that $\TM$ is, as an algebra, the full centre of the module category in the sense of Davydov, Kong and Runkel~\cite{MR3406516}, which for module categories of algebras is the full centre of Fjelstad, Fuchs, Runkel and Schweigert~\cite{MR2443266}. Consequently $\TM$ gives rise to Cardy algebras in the sense of Kong and Runkel~\cite{MR2551797}.
\end{abstract}

\maketitle

\section{Introduction}

Our primary aim in this paper is to establish a general, categorical relationship between the diagonal part of the modular invariant partition function and the corresponding Non-negative Integer Matrix representation (NIM-rep). The existence of such a relationship was first suggested by the appearance of an ADE pattern in the classification of such partition functions for $\mathfrak{su}(2)$ WZW-models~\cite{MR918402,MR1907190,Petkova2001} and subsequently by similar phenomena in the $\mathfrak{su}(n)$ case~\cite{MR1063590}. It should be noted that~\cite[Theorem 4.16, (2)]{Boeckenhauer1999ChiralSO} provides a prior description of this relationship, however, it is stated and proved within the operator algebra approach to conformal field theory (see also~\cite{MR1729094,MR1642584}).

To illustrate the phenomenon concretely, consider $\C_{\ell}$, the category of representations of the affine algebra $A_{1}^{(1)}$ at level $\ell$, which has $\ell + 1$ simple objects $x_0, \ldots, x_\ell$. The fusion algebra $F_{\ell}$ is the truncated $\mathrm{SU}(2)$ character ring, which is semisimple with spectrum
\begin{align*}
  \Spec(F_{\ell}) = \left\{ 2 \cos\left( \frac{\pi m}{h} \right) : m = 1, \ldots, h-1 \right\}
\end{align*}
where $h = \ell + 2$. Indecomposable module categories $\M_{\Q}$ over $\C_{\ell}$ are classified by ADE Dynkin graphs $\Q$ with Coxeter number $h$~\cite{MR2046203}. Let $\Nim_{\Q} \colon F_{\ell} \to \End(\Z^{\Q_{0}})$ denote the corresponding NIM-rep. The adjacency matrix of $\Q$ has eigenvalues
\begin{align*}
  \{ 2 \cos(\pi m / h) : m \in \mathrm{Exp}(\Q) \},
\end{align*}
where $\mathrm{Exp}(\Q)$ denotes the set of Coxeter exponents of $\Q$. Since the NIM-rep matrices $\Nim_{\Q}(x_{i})$ satisfy the same recurrence as the Chebyshev polynomials, the spectrum of $\Nim_{\Q}$ (as an $F_{\ell}$-module) is determined by the eigenvalues of the adjacency matrix, and one finds that
$\lambda_{m} = 2\cos(\pi m/h)$ appears in $\Nim_{\Q}$ with multiplicity equal to the multiplicity of $m$ as a Coxeter exponent.
On the other hand, the modular invariant partition functions for this category can also be related to the ADE Dynkin graphs~\cite{MR918402}: the diagonal entries $Z_{mm}$ of the corresponding modular invariant matrix record precisely the same Coxeter exponent multiplicities.
This is the relationship we seek to explain in general: that the diagonal part of the modular invariant is the spectrum of the NIM-rep (cf.~\cite[Section 4]{Gannon2005}).

We work within the framework established in~\cite{Hardiman:2019utc}, which starts with a modular tensor category $\C$ (which, from a physical point of view, describes the chiral data) and a module category $\M$ (which describes the boundary conditions), and then realises the partition function as a particular representation $\TM$ of the \emph{tube category}. The tube category, denoted $ \TC $, is a category whose objects coincide with those of $ \C $ and whose morphism are described by diagrams in $ \C $ inscribed upon the surface of a cylinder (see~\cite{hardiman_king} for a precise definition). Let $\RTC$ denote the category of representations of $\TC$, i.e.\ the category of contravariant functors from $\TC$ to the category of finite dimensional vector spaces. Crucially, the braiding on $\C$ provides a chain of equivalences going from $\RTC$, to the \emph{Drinfeld centre} $Z(\C)$, to the Deligne tensor product $\CC$, automatically labelling the simple objects in $\RTC$ by \emph{pairs} of simple objects in $\C$. This allows us to interpret representations of the tube category (i.e.\ objects in $\RTC$) as potential modular invariants. Indeed, the matrix form may be recovered by evaluating the representation's irreducible multiplicity spaces (see~\cite{hardiman19a} for a detailed description of this process).

We briefly establish some notation and conventions. We work over an algebraically closed field $\fld$ of characteristic zero. We write $\Sharp$ for the Yoneda embedding (which sends an object to its representable functor), so that for an object $X$ in $\TC$, $X\Sharp$ denotes the corresponding representable functor in $\RTC$. More generally, for an object/idempotent pair $(X,e)$ in the Karoubi envelope of $\TC$, $(X,e)\Sharp$ denotes the corresponding object in $\RTC$. In particular, the complete set of simple objects in $\RTC$ labelled by \emph{pairs} of simple objects in $\C$ (mentioned in the preceding paragraph) is given by
\begin{align*}
  \left\{\left(IJ,\e_{I}^{J}\right)\Sharp\right\}_{I,J \in \I}
\end{align*}
where $\e_{I}^{J}$ is the idempotent in $\TC$ specified in Definition~\ref{def:primitive-idem-tc}.
In the other direction, for an object $\F$ in $\RTC$, $\F\Flat$ denotes its underlying object in $\C$ (see Definition~\ref{def:flat}).

One advantage of establishing the partition function as a representation of the tube category is that the \say{diagonal part} of any object in $\RTC$ is canonically a representation of the fusion algebra $F$. To understand how this comes about, let $\tid$ denote the tensor identity in $\C$ and note that from \cite[Proposition 3.1.12]{hardiman2019module}, we have
\begin{align*}
  \Hom_{\RTC}\left( \tid\Sharp, \left(IJ,\e_{I}^{J}\right)\Sharp \right) = \Hom_{\TC}\left(\tid, \e_{I}^{J}\right) \cong \Hom_{\C}(\tid, IJ) \cong \delta_{J,I^{\vee}} \,\fld
\end{align*}
for all $I,J \in \Irr(\C)$, implying the following decomposition of $\tid\Sharp$, \emph{as an object in $\RTC$}:
\begin{align} \label{eq:tid_decomposition}
  \tid\Sharp = \bigoplus_{I \in \Irr(\C)} \left(II^{\vee},\e_{I}^{I^{\vee}}\right)\Sharp.
\end{align}

Furthermore, the endomorphism algebra $\End_{\TC}(\tid)$ is canonically isomorphic to the fusion algebra $F$ (see Remark~\ref{rem:isom_kc_endtid}). Therefore every representation $\F$ in $\RTC$ automatically defines an $F$-module, namely $\F(\tid)$ and, in light of~\eqref{eq:tid_decomposition}, this $F$-module may be thought of as the \say{diagonal part} of $\F$. The goal is therefore to identify this $F$-module with the NIM-rep when $\F = \TM$, and moreover to show that the $\lambda_{I}$-eigenspaces of $\TM(\tid)$ coincide with the diagonal multiplicity spaces $\TM_{I}^{I^{\vee}}$, thereby equating the entries $Z_{I,I^{\vee}}$ of the modular invariant with the NIM-rep multiplicities.

A word on conventions. With $Z_{IJ} = \dim \TM_{I}^{J}$ (see Section~\ref{sec:intro-tm}), the regular module category $\M = \C$ gives the charge-conjugation invariant $Z_{IJ} = \delta_{J,I^{\vee}}$ (Remark~\ref{rem:regular-module}), as in~\cite{MR1940282}, where the trivial algebra $A = \tid$ likewise gives charge conjugation. Accordingly, the entries $Z_{I,I^{\vee}}$ are the diagonal entries of $ZC$, where $C$ is the charge-conjugation matrix. They become the diagonal entries of $Z$ after the relabelling $J \mapsto J^{\vee}$, i.e.\ in conventions where the regular module category gives the identity matrix. This should be taken into account when comparing with~\cite{Boeckenhauer1999ChiralSO}. All simple objects of $\C_{\ell}$ are self-dual, so in the example above the two readings agree.

This ends up following from Theorem~\ref{thm:encircling-module-iso-nim-rep}, which proves a statement in a more general set-up, namely that the \emph{encircling module}, $\Enc$, is isomorphic to the NIM-rep, $\Nim$ (see Section~\ref{sec:mod-cat-fusion-alg} for definitions). Recall that $\Nim$ is defined for an arbitrary finite module category $\M \colon \C \to \End(\B)$ where $\C$ is a fusion category. In the case when $\C$ is equipped with a spherical pivotal structure, one may also define an additional $F$-module we call the encircling module, so named due the appearance of its action when depicted using the graphical calculus of monoidal categories~\eqref{eq:enc_mod}. Theorem~\ref{thm:encircling-module-iso-nim-rep} then establishes that, under the stipulation that $\M$ be pivotal (see Definition~\ref{def:pivotal-module-category}), $\Enc$ and $\Nim$ are isomorphic $F$-modules. It is fairly immediate that, under the additional assumptions required by the definition of $\TM$, the $F$-module arising from $\TM$ (the  \say{diagonal part}) is precisely $\Enc$ (see Remark~\ref{rem:encir-rep-mod-inv-diag}), while the identification of the primitive idempotents in $\End_{\TC}(\tid)$ with the spectral idempotents of $F$ (Proposition~\ref{prop:idempotents-agree}) shows that the $\lambda_{I}$-eigenspace of $\TM(\tid)$ is precisely $\TM_{I}^{I^{\vee}}$. The desired result, Corollary~\ref{cor:diag_tm_nim_rep}, therefore follows.

As a further application of the techniques developed in the preceding sections, we are also able to compute the quantum dimension $d(\TM)$. In~\cite{Hardiman:2019utc}, it was shown that $\TM$ is a modular invariant algebra provided the dimension condition $d(\TM) = d(\C)$ holds. However, this condition was required as an assumption. By expressing $d(\TM)$ in terms of the NIM-rep character $\chi_{\M}$ and applying the orthogonality of $\Spec(F)$, we show that $d(\TM) = \dim \TM_{\tid}^{\tid} \cdot d(\C)$ (Proposition~\ref{prop:dimension-tm}). Since indecomposability of $\M$ is equivalent to $\TM_{\tid}^{\tid} = \fld$, the dimension condition is automatically satisfied, and the assumption in~\cite{Hardiman:2019utc} may be removed (Corollary~\ref{cor:TM-mod-inv}).

Finally, to connect the $\TM$ construction to the existing literature, we show that under the equivalence $\RTC \simeq \cZ(\C)$, $\TM$ is isomorphic as an algebra to the \emph{full centre} $Z(\M)$~\cite{MR2443266,MR2658689,MR3406516} (Theorem~\ref{thm:tm-iso-full-centre}). In particular, for indecomposable $\M$, $\TM$ gives rise to a Cardy $\C$-algebra (Corollary~\ref{cor:tm-cardy-alg}). For this we also show that pivotality of $\M$ is equivalent to the existence of a module trace in the sense of Schaumann~\cite{MR3019263} (Proposition~\ref{prop:module-trace}).

In Section~\ref{sec:mod-cat-fusion-alg}, we introduce the encircling module $\Enc$ and prove Theorem~\ref{thm:encircling-module-iso-nim-rep}, which gives an isomorphism of $F$-modules between $\Enc$ and the NIM-rep $\Nim$. In Sections 3 and 4, we recall the tube category $\TC$ and the relationship between $\Spec(F)$ and the primitive idempotents in $\End_{\TC}(\tid)$. In Section 5, we combine these ingredients to establish our main result Corollary~\ref{cor:diag_tm_nim_rep}. In Section 6, we compute the dimension of $\TM$ (Proposition~\ref{prop:dimension-tm}), removing the dimension assumption from~\cite{Hardiman:2019utc}. Finally, in Section~\ref{sec:cft}, we identify $\TM$ with the full centre construction~\cite{MR2443266,MR2658689,MR3406516}.

\section{Modules over the Fusion Algebra}
\label{sec:mod-cat-fusion-alg}

Let $\fld$ be an algebraically closed field of characteristic zero. Let $\C$ be a fusion category over $\fld$ and let $F$ denote the ``complexified'' Grothendieck ring of $\C$; we call $F$ the \emph{fusion algebra} of $\C$ (cf.~\cite{Verlinde1988}). Let $ \M\colon \C \to \End(\B) $ be a finite\footnote{In the sense that $\B$ is semisimple with finitely many simple objects.} module category over $ \C $. Applying the functor $\Gr$, which takes a monoidal category to its corresponding Grothendieck ring, leads to the $\Gr(\C)$-module
\begin{align*}
  \Gr(\C) \to \Gr(\End(\B)) = \End(\Gr(\B))
\end{align*}
known in the literature as the NIM-rep. One may also consider the ``complexification'' of this module, i.e.\ the $F$-module
\begin{align*}
  \Nim = \Gr(\B) \otimes \fld.
\end{align*}

Suppose in addition that $\C$ is equipped with a spherical pivotal structure given by natural isomorphisms $p_{X} \colon X \to X^{\vee \vee}$. Then $\M$ gives rise to another $F$-module as follows. Let $\Enc$ be the vector space $\Nat(\id_{\B})$ equipped with the $F$-action
\begin{align*}
  X \cdot \alpha = \M(\ann_{X^{\vee}}) \circ \M(p_{X}) \otimes \alpha \otimes \id_{\M(X)^{\vee}} \circ \M(\cre_{X}).
\end{align*}
We can draw this action with respect to the graphical calculus of $\End(\B)$:
\begin{align} \label{eq:enc_mod}
  X \cdot \alpha = \hspace{-0.6em}
  \begin{array}{c}
    \begin{tikzpicture}[scale=0.15,every node/.style={inner sep=0,outer sep=-1}]
      \node [blue] at (-5,3) {$ X $};
      \draw [thick, blue] (0,0) ellipse (4 and 4);
      \draw[fill=blue] (-4,0) circle (20pt);
      \node [draw,outer sep=0,inner sep=1,minimum height=12,minimum width=24,fill=white] (v5) at (0,0) {$  \alpha $};
    \end{tikzpicture}
  \end{array}
\end{align}
where the colour blue indicates evaluation under $\M$ and the image of the pivotal structure on $\C$ is depicted by a blue dot. In accordance with this figure, we call $\Enc$ the \emph{encircling module} over $F$. Identifying the encircling module $\Enc$ with the NIM-rep $\Nim$ requires an additional assumption on $\M$, which will now be introduced via the following two definitions.

\begin{DEF}
  \label{def:full-sub-cat}
  Let $\D$ be the idempotent completion of the image of $\M$, that is, the full subcategory of $\End(\B)$ consisting of subobjects of objects of the form $\M(X)$ for $X \in \C$.
\end{DEF}

Throughout, the monoidal product of $\End(\B)$ is composition of functors, $G \otimes H = G \circ H$, so that $\M(X \otimes Y) \cong \M(X) \circ \M(Y)$. For $i,j \in \Irr(\B)$ let $\simp{i}{j}$ be the simple object in $\End(\B)$ that takes $j$ to $i$ and $k$ to $0$ for $k \neq j$. Thus $\Hom_{\End(\B)}(\simp{i}{j},G) \cong \Hom_{\B}(i,G(j))$ for $G \in \End(\B)$, and $\simp{i}{j} \otimes \simp{j}{k} = \simp{i}{k}$. We note that $ \{\simp{i}{j}\}$ forms a complete set of simple objects in $\End(\B)$ and that the object $\simp{j}{i}$ canonically admits the structure of a left and right dual to $\simp{i}{j}$. Note that $\simp{i}{j} \in \D$ if and only if there exists $X \in \C$ such that $$\Hom_{\End(\B)}(\simp{i}{j},\M(X)) \neq 0.$$
Furthermore, such $\simp{i}{j}$ form a complete set of simple objects in $\D$. Let $\J$ be the set of pairs $(i,j)$ such that $\simp{i}{j}$ is in $\D$. As $\M$ is monoidal, for all $i,j \in \Irr(\B)$, $\J$ satisfies the following properties:
\begin{align}\hspace{-4.5em} \label{eq:j-props}
  \left\{\mystrut(27,0)\right.\hspace{-3.5em}
  \begin{tabular}{p{.8\textwidth}}
    \begin{itemize}
    \item[] $(i,i) \in \J$, \hspace{10.75em} as $e_{ii} < \tid_{\D} = \M(\tid_{\C}),$
    \item[] $(i,j) \in \J \implies (j,i) \in \J$, \hspace{3.75em} as $e_{ji}$ is dual to $e_{ij},$
    \item[] $(i,j), (j,k) \in \J \implies (i,k) \in \J$, \hspace{1em} as $e_{ij} \otimes e_{jk} = e_{ik}.$
    \end{itemize}
  \end{tabular}
\end{align}

\begin{DEF}
  \label{def:pivotal-module-category}
  A module category $\M \colon \C \to \End(\B)$ is called \emph{pivotal} if the image of the pivotal structure on $\C$ induces a pivotal structure on $\D$.
\end{DEF}

To clarify this definition, observe that, by functoriality, $ \M(p_{X}) $ is automatically natural with respect to morphisms in $ \C $. However, to give a pivotal structure on $ \D $ it must be natural with respect \emph{all} morphisms in $ \D $. In other words, the diagram
\begin{equation*}
  \begin{tikzcd}[column sep=3.8em]
    \M(X) \arrow{r}{\alpha} \arrow[swap]{d}{\M(p_{X}) \ }&
    \M(Y) \arrow{d}{\ \M(p_{Y})}\\
    \M(X^{\vee\vee}) \arrow{r}{\alpha^{\vee\vee}}&
    \M(Y^{\vee\vee})
  \end{tikzcd}
\end{equation*}
must commute for all $ \alpha \in \Hom_{\D}(\M(X),\M(Y)) $. Note that $\M(X^{\vee})$ is canonically dual to $\M(X)$, i.e.\ canonically isomorphic to $\M(X)^{\vee}$.

\begin{EX}
  \cite[Section 8]{Hardiman:2019utc} shows that every module category over the fixed level representations of the affine algebra $A_{1}^{(1)}$ is pivotal.
\end{EX}

\begin{EX}
  \label{ex:non-pivotal}
  Let $ \C $ be the fusion category of representations of $ \mathbb{Z}/2\mathbb{Z} $. A complete set of simples in $ \C $ is given by $ \{\underline{0}, \underline{1} \} $, where $ \underline{0} $ is the tensor unit and $ \underline{1} \otimes \underline{1} = \underline{0} $. We may then equip $ \C $ with the pivotal structure $ \delta_{\underline{1}} = - \id_{\underline{1}} $ (so that $ d(\underline{1}) = -1 $). One may then check that the module category
  \begin{align*}
    \M\colon \C &\to \Vect \\
    \underline{0} &\mapsto \mathbb{K} \\
    \underline{1} &\mapsto \mathbb{K},
  \end{align*}
  fails to be pivotal.
\end{EX}

Now that the condition that a module category be pivotal has been defined, we may proceed with the proof that equates the encircling module $\Enc$ with the NIM-rep $\Nim$.

Let $p_{A},q_{A}\colon A \to A^{\vee\vee}$ denote two pivotal structures on $\D$. To compare $p_{A}$ and $q_{A}$ it is useful to consider the following numbers,
\begin{align*}
  \mu_{ij} = p_{\simp{i}{j}}^{-1} \circ q_{\simp{i}{j}}^{\vphantom{-1}} \in \Aut_{\D}(\simp{i}{j}) = \fld^{*}, \quad \text{for $(i,j) \in \J$.}
\end{align*}
Appendix~\ref{app:compare-pivotal-struct} establishes certain properties that these numbers satisfy (within the more general context of an arbitrary multifusion category). In particular, we have
\begin{enumerate}[label=(\roman*)]
\item $\mu_{ii}=1$ \label{it:ii}
\item $\mu_{ij} = \mu_{ji}^{-1}$ \label{it:ij}
\item $\mu_{ij}\mu_{jk} = \mu_{ik}.$ \label{it:ijjk}
\end{enumerate}
where the existence of the relevant numbers is guaranteed by~\eqref{eq:j-props}.
\begin{LEMMA}
  \label{lem:compare-pivotal-structs}
  There exists a set of scalars $\{\lambda_{i}\}_{i \in \Irr(\B)}$ such that
  \begin{align} \label{eq:mu-equation}
    \mu_{ij} = \frac{\lambda_{i}}{\lambda_{j}}.
  \end{align}
  \proof
  By~\eqref{eq:j-props}, $\J$ is an equivalence relation on $\Irr(\B)$. Choose a representative $i_{0}$ of each equivalence class, and set $\lambda_{i} = \mu_{i i_{0}}$ for $i$ in the class of $i_{0}$. If $(i,j) \in \J$, then $i$ and $j$ lie in the same class, with representative $i_{0}$ say, and by~\ref{it:ij} and~\ref{it:ijjk},
  \begin{align*}
    \frac{\lambda_{i}}{\lambda_{j}} = \mu_{i i_{0}} \mu_{i_{0} j} = \mu_{ij}.
  \end{align*}
  (Equivalently, by~\ref{it:ii}--\ref{it:ijjk}, $\mu$ is a $\fld^{*}$-valued $1$-cocycle on the clique complex of the graph defined by $\J$. Each component of this complex is a simplex, so its first cohomology vanishes.)
  \endproof
\end{LEMMA}

\begin{REM}
  \label{rem:multiple-solutions}
  The solution of~\eqref{eq:mu-equation} is unique up to rescaling on each equivalence class of $\J$. Indeed, if $\{\lambda_i\}$ and $\{\lambda'_i\}$ are two solutions, then $\lambda'_i/\lambda_i = \lambda'_j/\lambda_j$ whenever $(i,j) \in \J$. Hence $\lambda'_i = c_{[i]} \lambda_i$ for scalars $c_{[i]} \in \fld^{*}$ depending only on the class $[i]$ of $i$. The classes of $\J$ are precisely the sets of simple objects of the indecomposable summands of $\M$. So the solution is unique up to a single global scalar if and only if $\M$ is indecomposable.
\end{REM}

\begin{THM}
  \label{thm:encircling-module-iso-nim-rep}
  Let $\M\colon \C \to \End(\B)$ be a finite pivotal module category. Then $\Nim$ and $\Enc$ are isomorphic $F$-modules.
  \proof
  Let $\{\lambda_{i}\}$ be a set of scalars satisfying the statement of Lemma~\ref{lem:compare-pivotal-structs} in the case when $p_{A}$ is the pivotal structure induced by $\M$ and $q_{A}$ is the standard pivotal structure on $\End(\B)$. We claim that the following map is an isomorphism of $F$-modules,
  \begin{align*}
    \Phi  \colon \Enc &\to \Nim\\
    \phi^{i} &\mapsto \lambda_{i} \ i
  \end{align*}
  where $\phi^{i}$ is defined by $\phi^{i}_{j} = \begin{cases} \id_{i} & \text{if } i=j \\ 0 & \text{else}\end{cases} \quad \forall j \in \Irr(\B)$.

  As elements of the form $\phi^{i}$ form a basis of $\Nat(\id_{\B})$, $\Phi$ is clearly an isomorphism of vector spaces. Therefore, we only have to check that, for all $X \in F$, the diagram
  \begin{equation*}
    \begin{tikzcd}[column sep=3.8em]
      \Nat(\id_{\B}) \arrow{r}{X \cdot \Blank} \arrow[swap]{d}{\Phi \ }&
      \Nat(\id_{\B})  \\
      \Gr(\B) \arrow{r}{X \cdot \Blank}&\arrow[swap]{u}{\ \Phi^{-1}}
      \Gr(\B)
    \end{tikzcd}
  \end{equation*}
  commutes, i.e. that for all $i \in \Irr(\B)$, 
  \begin{align*}
    X \cdot \phi^{i} &= \lambda_{i} \ \Phi^{-1}\big(\M(X)(i) \big)\\
                     &= \lambda_{i} \sum_{j \in \Irr(\B)} \my-hom\big(j,\M(X)(i)\big) \Phi^{-1}\big(j \big) \\
                     &= \sum_{j\in \Irr(\B)} \frac{\lambda_{i}}{\lambda_{j}} \my-hom\big(j,\M(X)(i)\big) \phi^{j}\\
                     &= \sum_{j\in \Irr(\B)} \frac{\lambda_{i}}{\lambda_{j}} \my-hom\big(\simp{j}{i}, \M(X)\big) \phi^{j},
  \end{align*}
  where $\my-hom$ is used to denote the dimension of the corresponding $\Hom$-space.
  We compute,
  \begin{align*}
    X \cdot \phi^{i} =
    \hspace{-0.6em}
    \begin{array}{c}
      \begin{tikzpicture}[scale=0.15,every node/.style={inner sep=0,outer sep=-1}]
        \node [blue] at (-5,3) {$ X $};
        \draw [thick, blue] (0,0) ellipse (4 and 4);
        \node [draw,outer sep=0,inner sep=1,minimum height=12,minimum width=14,fill=white] at (-4,0) {$p$};
        \node [draw,outer sep=0,inner sep=1,minimum height=12,minimum width=14,fill=white] (v5) at (0,0) {$  \phi^{i} $};
      \end{tikzpicture}
    \end{array}
    &= \sum_{j \in \Irr(\B)} \hspace{-0.6em}
      \begin{array}{c}
        \begin{tikzpicture}[scale=0.15,every node/.style={inner sep=0,outer sep=-1}]
          \node [blue] at (-5,3) {$ X $};
          \draw [thick, blue] (0,0) ellipse (4 and 4);
          \node [draw,outer sep=0,inner sep=1,minimum height=12,minimum width=14,fill=white] at (-4,0) {$p$};
          \node [draw,outer sep=0,inner sep=1,minimum height=12,minimum width=14,fill=white] (v5) at (0,0) {$  \phi^{i} $};
          \node [draw,outer sep=0,inner sep=1,minimum height=12,minimum width=14,fill=white] (v5) at (-8,0) {$  \phi^{j} $};
        \end{tikzpicture}
      \end{array}
    \\
    &= \sum_{j \in \Irr(\B)} \my-hom\big(\simp{j}{i}, \M(X)\big) \hspace{-0em}
      \begin{array}{c}
        \begin{tikzpicture}[scale=0.15,every node/.style={inner sep=0,outer sep=-1}]
          \node [blue] at (-4.7,3.7) {$ \simp{j}{i} $};
          \draw [thick, blue] (0,0) ellipse (4 and 4);
          \node [draw,outer sep=0,inner sep=1,minimum height=14,minimum width=16,fill=white] at (-4,0) {$p_{\simp{j}{i}}$};
        \end{tikzpicture}
      \end{array}
    \\
    &= \sum_{j\in \Irr(\B)} \frac{\lambda_{i}}{\lambda_{j}} \my-hom\big(\simp{j}{i}, \M(X)\big) 
      \begin{array}{c}
        \begin{tikzpicture}[scale=0.15,every node/.style={inner sep=0,outer sep=-1}]
          \node [blue] at (-4.7,3.7) {$ \simp{j}{i} $};
          \draw [thick, blue] (0,0) ellipse (4 and 4);
          \node [draw,outer sep=0,inner sep=1,minimum height=14,minimum width=16,fill=white] at (-4,0) {$q_{\simp{j}{i}}$};
        \end{tikzpicture}
      \end{array}
    \\
    &= \sum_{j\in \Irr(\B)} \frac{\lambda_{i}}{\lambda_{j}}  \my-hom\big(\simp{j}{i}, \M(X)\big) \phi^{j}
  \end{align*}
  where the penultimate equality is due to Lemma~\ref{lem:compare-pivotal-structs} and the final equality is due to the fact that, using the standard pivotal structure on $\End(\B)$, $d(\simp{j}{i})=1$. \endproof
\end{THM}

A special case of Theorem~\ref{thm:encircling-module-iso-nim-rep} was proved by~\cite[Proposition 8.5]{Hardiman:2019utc}. In this prior proof, the $\lambda_{i}$ numbers are defined, not in terms of the image of the pivotal structure on $\C$, but rather as the entries of the Perron-Frobenius eigenvector of the quiver used to construct $\M$. The following proposition shows that, in general, the $\lambda_{i}$ form an eigenvector of the NIM-rep for the character $d$. Remark~\ref{rem:perron-frobenius} explains when this identifies them with a Perron--Frobenius eigenvector.

\begin{PROP}
  \label{prop:identify-numbers}
  Let $\M \colon \C \to \End(\B)$ be a finite pivotal module category and let $\{\lambda_{i}\}$ be a set of scalars as in the proof of  Theorem~\ref{thm:encircling-module-iso-nim-rep}. Then the vector $$v=\left(\lambda_{i}\right)_{i \in \Irr(\B)} \in \Gr(\B),$$ lies in the $d$-eigenspace of $\Nim$ (where $d \in \Spec(F)$ is the categorical trace on $\C$, see Definition~\ref{def:spherical-cat}).

  \proof Let $\Phi  \colon \Nat(\id_{\B}) \to \Gr(\B)$ be the isomorphism from $\Enc$ to $\Nim$ described within the proof of Theorem~\ref{thm:encircling-module-iso-nim-rep} and let $\underline{\id}$ denote the identity element in $\Nat(\id_{\B})$. By the definition of $\Enc$, we have, for $X \in F$,
  \begin{align*}
    X \cdot \underline{\id} = 
    \begin{array}{c}
      \begin{tikzpicture}[scale=0.1,every node/.style={inner sep=0,outer sep=-1}]
        \node [blue] at (-5.5,3) {$ X $};
        \draw [thick, blue] (0,0) ellipse (4 and 4);
        \draw[fill=blue] (-4,0) circle (20pt);
      \end{tikzpicture}
    \end{array}
    = d(X) \ \underline{\id},
  \end{align*}
  implying that $\underline{\id}$ lies within the $d$-eigenspace of $\Enc$. Furthermore, as
  \begin{align*}
    \Phi(\underline{\id}) = \Phi\left(\sum_{i} \phi^{i} \right) = \sum_{i}  \Phi\left( \phi^{i} \right) = \sum_{i} \lambda_{i} \ i = v,
  \end{align*}
  we deduce that $v$ lies within the $d$-eigenspace of $\Nim$. \endproof
\end{PROP}

\begin{REM}
  \label{rem:perron-frobenius}
  A $d$-eigenvector of $\Nim$ need not be a Perron--Frobenius eigenvector. To identify $v$ with the Perron--Frobenius eigenvector one needs, in addition, that the $d$-eigenspace of $\Nim$ is one-dimensional and that $d$ is the Frobenius--Perron dimension. Suppose that $\fld = \mathbb{C}$ and $d(X) = \mathrm{FPdim}(X)$ for all $X \in \C$ (for instance, $\C$ pseudo-unitary with its canonical spherical structure~\cite{MR2183279}, as is the case for $\C_{\ell}$). Suppose also that $\M$ is indecomposable. Then the Frobenius--Perron vector $w$ of $\Nim$, which has positive entries and satisfies $X \cdot w = \mathrm{FPdim}(X)\, w$ for all $X$~\cite[Chapter~3]{Etingof15}, lies in the $d$-eigenspace. Since $\M$ is indecomposable, every entry of the matrix $P = \sum_{S \in \I} \Nim(S)$ is positive. By the Perron--Frobenius theorem, the eigenspace of $P$ for the eigenvalue $\sum_{S} d(S)$ of the positive eigenvector $w$ is therefore one-dimensional. This eigenspace contains the $d$-eigenspace of $\Nim$, so $v$ is a scalar multiple of $w$. This recovers the description of the $\lambda_{i}$ in~\cite[Proposition 8.5]{Hardiman:2019utc}. When $\C$ is modular, the one-dimensionality of the $d$-eigenspace for indecomposable $\M$ also follows from Proposition~\ref{prop:dimension-tm}.
\end{REM}

The scalars $\lambda_{i}$ have a further interpretation, in terms of the \emph{module traces} of Schaumann~\cite{MR3019263}. Recall that a trace $\Theta$ on $\B$ is a family of linear maps $\Theta_{m} \colon \End_{\B}(m) \to \fld$ that is symmetric and non-degenerate. It is a \emph{module trace} if moreover $\Theta_{X \otimes m} = \Theta_{m} \circ \tr_{X}$ for all $X \in \C$ and $m \in \B$, where $\tr_{X} \colon \End_{\B}(X \otimes m) \to \End_{\B}(m)$ is the partial trace formed using the pivotal structure of $\C$~\cite[(3.24) and Definition~3.7]{MR3019263}. (The definitions of~\cite{MR3019263} are stated over $\mathbb{C}$, but make sense over $\fld$.) As $\B$ is semisimple, every trace on $\B$ is of the form
\begin{align*}
  \Theta^{\nu}_{m}(f) = \sum_{i \in \Irr(\B)} \nu_{i} \, \tr\big(\Hom_{\B}(i,f)\big)
\end{align*}
for a unique $\nu \colon \Irr(\B) \to \fld^{*}$~\cite[Proposition~5.2]{MR3019263}. Let $q$ denote the standard pivotal structure on $\End(\B)$. For $\nu$ as above, let $q^{\nu}$ be the pivotal structure on $\D$ obtained by composing $q$ with the monoidal automorphism of $\id_{\D}$ which acts on $\simp{i}{j}$ by $\nu_{j}/\nu_{i}$, so that $q^{\nu}_{\simp{i}{j}} = \frac{\nu_{j}}{\nu_{i}} q_{\simp{i}{j}}$.

\begin{PROP}
  \label{prop:module-trace}
  Let $\M \colon \C \to \End(\B)$ be a finite module category over a spherical fusion category $\C$, and let $\nu \colon \Irr(\B) \to \fld^{*}$. Then $\Theta^{\nu}$ is a module trace if and only if $\M(p_{X}) = q^{\nu}_{\M(X)}$ for all $X \in \C$. In particular, $\M$ is pivotal if and only if $\B$ admits a module trace. In this case the dimensions $\nu_{i} = \Theta^{\nu}_{i}(\id_{i})$ of the simple objects form a solution $\{\lambda_{i}\}$ of~\eqref{eq:mu-equation}.
  \proof
  Since the $\nu_{i}$ are non-zero, $\Theta^{\nu}$ is symmetric and non-degenerate~\cite[Lemma~5.1]{MR3019263}, so only the compatibility condition $\Theta^{\nu}_{X \otimes m} = \Theta^{\nu}_{m} \circ \tr_{X}$ is at issue. By additivity it suffices to consider $m = j \in \Irr(\B)$. We have
  \begin{align*}
    \End_{\B}(\M(X)(j)) = \End_{\End(\B)}(\M(X) \simp{j}{j}), \qquad \M(X)\simp{j}{j} \cong \bigoplus_{i} V_{ij} \otimes \simp{i}{j},
  \end{align*}
  where $V_{ij} = \Hom_{\End(\B)}(\simp{i}{j},\M(X))$. So an endomorphism $f$ of $\M(X)(j)$ is a family of linear maps $f_{i} \in \End(V_{ij})$, and $\Theta^{\nu}_{X \otimes j}(f) = \sum_{i} \nu_{i} \tr(f_{i})$.

  Let $\theta_{X} = q_{\M(X)}^{-1} \circ \M(p_{X}) \in \Aut(\M(X))$, and let $\theta_{X}^{ij} \in \Aut(V_{ij})$ be its action on $V_{ij}$. The partial trace $\tr_{X}(f)$ is the left partial trace of $f$ over $\M(X)$ in $\End(\B)$, formed using $\M(p_{X}) = q_{\M(X)} \circ \theta_{X}$. Every $\simp{i}{j}$ has dimension $1$ with respect to $q$, so
  \begin{align*}
    \tr_{X}(f) = \sum_{i} \tr\big(f_{i} \circ (\theta_{X}^{ij})^{-1}\big) \, \id_{j}.
  \end{align*}
  Hence $\Theta^{\nu}$ is a module trace if and only if $\sum_{i} \nu_{i} \tr(f_{i}) = \nu_{j} \sum_{i} \tr\big(f_{i} \circ (\theta^{ij}_{X})^{-1}\big)$ for all $X$, $j$ and $f_{i}$. Equivalently, $\theta_{X}^{ij} = \frac{\nu_{j}}{\nu_{i}} \id_{V_{ij}}$ for all $X$, $i$ and $j$, i.e.\ $\M(p_{X}) = q^{\nu}_{\M(X)}$ for all $X$.

  If this holds, then $q^{\nu}$ is a pivotal structure on $\D$ extending $\M(p)$, so $\M$ is pivotal. As a pivotal structure on $\D$ is determined by its values on the objects $\M(X)$, the pivotal structure induced by $\M$ is then $q^{\nu}$, and $\mu_{ij} = (q^{\nu}_{\simp{i}{j}})^{-1} \circ q_{\simp{i}{j}} = \nu_{i}/\nu_{j}$. Conversely, suppose that $\M$ is pivotal with induced pivotal structure $\tilde{p}$ on $\D$, and let $\{\lambda_{i}\}$ be as in Lemma~\ref{lem:compare-pivotal-structs} for $p_{A} = \tilde{p}$ and $q_{A} = q$. Then $\tilde{p}_{\simp{i}{j}} = \mu_{ij}^{-1} q_{\simp{i}{j}} = \frac{\lambda_{j}}{\lambda_{i}} q_{\simp{i}{j}}$, so $\tilde{p} = q^{\lambda}$ by naturality. Hence $\Theta^{\lambda}$ is a module trace. \endproof
\end{PROP}

For example, the module category of Example~\ref{ex:non-pivotal} admits no module trace, in agreement with~\cite[Example~3.13]{MR3019263}.

\section{Preliminaries on the tube category}

Let $ \fld $ be an algebraically closed field of characteristic zero and let $ \C $ be a spherical fusion category over $ \fld $ with complete set of simples $ \I $.

\begin{DEF} \label{def:spherical-cat}
  Recall that a \emph{spherical category} is a pivotal category in which the left and right categorical traces coincide. In particular, the spherical structure on $\C$ defines a \emph{categorical dimension} function $d \colon \K(\C) \to \fld$ by setting $d(X) = \tr(\id_X)$ for simple $X$ and extending linearly. The \emph{global dimension} of $\C$ is defined by $d(\C) = \sum_{I \in \I} d(I)^2$.
\end{DEF}

The \emph{tube category} of $ \C $, denoted $ \TC $, is a category whose objects coincide with those of $ \C $ and whose $ \Hom $-spaces are given by
\begin{align*}
  \Hom_{\TC}(X,Y) \defeq \bigoplus\limits_S \Hom_{\C}(SX,YS)
\end{align*}
where, as per our conventions, the direct sum ranges over $ \I $ and the monoidal product symbol is suppressed (see~\cite{Lan2025} for an alternative characterisation). To depict a morphism in $ \TC $ using the graphical calculus of spherical fusion categories we take $ \alpha \in \Hom_{\C}(GX,YG) $ and write
\begin{align} \label{eq:tc_morphism}
  \alpha_G =
  \begin{array}{c}
    \begin{tikzpicture}[scale=0.25,every node/.style={inner sep=0,outer sep=-1}]
      \node (v1) at (0,4) {};
      \node (v4) at (0,-4) {};
      \node (v2) at (4,0) {};
      \node (v3) at (-4,0) {};
      \node (v5) at (-2,2) {};
      \node (v6) at (2,-2) {};
      \node (v7) at (2,2) {};
      \node (v11) at (-2,-2) {};
      \node [draw,diamond,outer sep=0,inner sep=.5,minimum size=22,fill=white] (v9) at (0,0) {$ \alpha $};
      \node at (3,3) {$ X $};
      \node at (-3,-3) {$ Y $};
      \node at (-3,3) {$ G $};
      \node at (3,-3) {$ G $};
      \draw [thick] (v9) edge (v6);
      \draw [thick] (v9) edge (v7);
      \draw [thick] (v9) edge (v11);
      \draw [thick] (v5) edge (v9);
      \draw[very thick, red]  (v1) edge (v3);
      \draw[very thick, red]  (v2) edge (v4);
      \draw[very thick]  (v1) edge (v2);
      \draw[very thick]  (v3) edge (v4);
    \end{tikzpicture}
  \end{array}
\end{align}
as shorthand for $ \bigoplus_S\sum_b (\id_Y \otimes \ b^*) \circ \alpha \circ (b \otimes \id_X) \in \Hom_{\TC}(X,Y) $, where $ \{ b \} $ is a basis of $ \Hom_{\C}(S,G) $ and $ \{ b^* \} $ is the corresponding dual basis of $ \Hom_{\C}(G,S) $ with respect to the perfect pairing given by composition into $ \End_{\C}(S) = \fld $, see \cite[Proposition 3.1]{hardiman_king}. The intuition is that whereas morphisms in $ \C $ may be represented graphically as diagrams drawn on a bounded region of the plane, morphisms in $ \TC $ are given  by diagrams drawn on a cylinder. In particular, the red lines in~\eqref{eq:tc_morphism} should be thought of as being glued; this is compatible with our notation as one may indeed show that
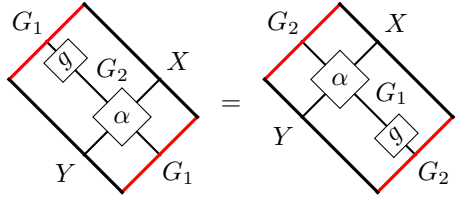
\begin{align*} 
  \begin{array}{c}
    \begin{tikzpicture}[scale=0.25,every node/.style={inner sep=0,outer sep=-1}]
      \node (v1) at (-2,6) {};
      \node (v4) at (0,-4) {};
      \node (v2) at (4,0) {};
      \node (v3) at (-6,2) {};
      \node (v5) at (-4,4) {};
      \node (v6) at (2,-2) {};
      \node (v7) at (2,2) {};
      \node (v11) at (-2,-2) {};
      \node [draw,diamond,outer sep=0,inner sep=.5,minimum size=22,fill=white] (v9) at (0,0) {$ \alpha $};
      \node at (3,3) {$ X $};
      \node at (-3,-3) {$ Y $};
      \node at (-5,5) {$ G_1 $};
      \node at (3,-3) {$ G_1 $};
      \node at (-0.5,2.5) {$ G_2 $};
      \draw [thick] (v9) edge (v6);
      \draw [thick] (v9) edge (v7);
      \draw [thick] (v9) edge (v11);
      \draw [thick] (v5) edge (v9);
      \node [draw,rotate=45,outer sep=0,inner sep=2,minimum width=13,fill=white] (v8) at (-3,3) {$ g $};
      \draw[very thick, red]  (v1) edge (v3);
      \draw[very thick, red]  (v2) edge (v4);
      \draw[very thick]  (v1) edge (v2);
      \draw[very thick]  (v3) edge (v4);
    \end{tikzpicture}
  \end{array}
  =
  \begin{array}{c}
    \begin{tikzpicture}[scale=0.25,every node/.style={inner sep=0,outer sep=-1}]
      \node (v1) at (-2,6) {};
      \node (v4) at (0,-4) {};
      \node (v2) at (4,0) {};
      \node (v3) at (-6,2) {};
      \node (v5) at (-4,4) {};
      \node (v6) at (2,-2) {};
      \node (v7) at (0,4) {};
      \node (v11) at (-4,0) {};
      \node [draw,diamond,outer sep=0,inner sep=.5,minimum size=22,fill=white] (v9) at (-2,2) {$ \alpha $};
      \node at (1,5) {$ X $};
      \node at (-5,-1) {$ Y $};
      \node at (-5,5) {$ G_2 $};
      \node at (3,-3) {$ G_2 $};
      \node at (0.75,1.25) {$ G_1 $};
      \draw [thick] (v9) edge (v6);
      \draw [thick] (v9) edge (v7);
      \draw [thick] (v9) edge (v11);
      \draw [thick] (v5) edge (v9);
      \node [draw,rotate=45,outer sep=0,inner sep=2,minimum width=13,fill=white] (v8) at (1,-1) {$ g $};
      \draw[very thick, red]  (v1) edge (v3);
      \draw[very thick, red]  (v2) edge (v4);
      \draw[very thick]  (v1) edge (v2);
      \draw[very thick]  (v3) edge (v4);
    \end{tikzpicture}
  \end{array}
\end{align*}
for any $ \alpha \in \Hom_{\C}(G_2X,YG_1) $ and $ g \in \Hom_{\C}(G_1,G_2) $. Composition in $ \TC $ is then defined following the intuition of vertically stacking the cylinders:
\begin{align} \label{eq:composition_in_tc}
  \beta_H \circ \alpha_G \defeq \bigoplus\limits_{T} \sum\limits_{b}
  \begin{array}{c}
    \begin{tikzpicture}[scale=0.5,every node/.style={inner sep=0,outer sep=-1}]
      \node (v1) at (-1.25,2.75) {};
      \node (v4) at (-1.25,-5.25) {};
      \node (v2) at (2.75,-1.25) {};
      \node (v3) at (-5.25,-1.25) {};
      \node (v9) at (0.75,0.75) {};
      \node (v6) at (-3.25,-3.25) {};
      \node (v70) at (0.75,-3.25) {};
      \node (v7) at (-3.25,0.75) {};
      \node at (0.75,-1.25) {$ G $};
      \node at (-3.25,-1) {$ H $};
      \node at (1.25,1.25) {$ X $};
      \node at (-3.75,-3.75) {$ Z $};
      \node at (-1.5607,-1.0119) {$ Y $};
      \node at (1.2,-3.7) {$ T $};
      \node at (-3.75,1.25) {$ T $};
      \node (v11) at (-2.625,-1.875) {};
      \node (v12) at (-1.875,-2.625) {};
      \draw [thick] (v11) edge (v12);
      \node (v14) at (-0.625,0.125) {};
      \node (v13) at (0.125,-0.625) {};
      \draw [thick] (v13) edge (v14);
      \node (v15) at (-2.45,-0.425) {};
      \node (v16) at (-2.075,-0.05) {};
      \draw [thick] (v15) to[out=-45,in=135] (v11);
      \draw [thick] (v16) to[out=-45,in=135] (v14);
      \node (v17) at (-0.05,-2.075) {};
      \node (v18) at (-0.425,-2.45) {};
      \draw [thick] (v13) to[out=-45,in=135] (v17);
      \draw [thick] (v12) to[out=-45,in=135] (v18);
      \node [diamond,draw,outer sep=0,inner sep=0.3,minimum size=25,fill=white] (v50) at (-2.25,-2.25) {\mbox{$ \beta $}};
      \node [diamond,draw,outer sep=0,inner sep=-0.2,minimum size=25,fill=white] (v5) at (-0.25,-0.25) {\mbox{$ \alpha $}};
      \node [draw,rotate=45,outer sep=0,inner sep=2,minimum height=10,minimum width=13,fill=white] (v8) at (-2.5,0) {$ b $};
      \node [draw,rotate=45,outer sep=0,inner sep=2,minimum height=10,minimum width=13,fill=white] (v10) at (0,-2.5) {$ b^* $};
      \draw[thick]  (v50) edge (v5);
      \draw [thick] (v50) edge (v6);
      \draw[thick]  (v7) edge (v8);
      \draw[thick]  (v10) edge (v70);
      \draw[thick]  (v9) edge (v5);
      \draw[very thick, red]  (v1) edge (v3);
      \draw[very thick, red]  (v2) edge (v4);
      \draw[very thick]  (v1) edge (v2);
      \draw[very thick]  (v3) edge (v4);
    \end{tikzpicture}
  \end{array}.
\end{align}
This intuition, together with the associativity of the tensor product, guarantee that composition in $ \TC $ is associative.

\begin{REM} \label{rem:isom_kc_endtid}
  Recall that $F$ denotes the ``complexified'' Grothendieck ring of $\C$, a.k.a. the fusion algebra. $ \End_{\TC}(\tid) $ and $ F $ are canonically isomorphic algebras. Indeed, $ \End_{\TC}(\tid) = \bigoplus_S \End(S) = \bigoplus_S \fld $ is precisely the underlying vector space of $ F $. Furthermore, composition in $ \End_{\TC}(\tid) $ corresponds to the tensor product in $ F $.
\end{REM}

\begin{REM} \label{rem:equiv_zc}
  As $ \C $ is a fusion category the Yoneda embedding gives an equivalence between $ \C $ and $ \RC $. As described in~\cite[Section 7]{hardiman19a}, the data required to extend the image of $ X $ under the Yoneda embedding to $ \TC $ corresponds to a half braiding on $ X $. Combining these facts yields an equivalence between $ Z(\C) $  and $ \RTC $, where $ Z(\C) $ is the \emph{Drinfeld centre} of $ \C $ (see also~\cite{Mueger2003,MR1990929,KnoetzeleST2024}).
\end{REM}

\section{$\Spec F$ and the tube category}

We now take $\C$ to be a \emph{modular tensor category} i.e.\ $\C$ is a spherical fusion category (as in the previous section) equipped with a balanced, non-degenerate braiding. As before, let $F$ denote the fusion algebra of $\C$ and let $\I$ denote a complete set of simple objects in $\C$. Recall that the (unnormalised) $S$-matrix satisfies $\Sc^{2} = d(\C)\, C$, where $C_{IJ} = \delta_{J,I^{\vee}}$ is the charge-conjugation matrix~\cite[Section 8.14]{Etingof15}; as $\Sc$ is invertible, in particular $d(\C) \neq 0$. We recall that the \emph{spectrum} of $ F $, $ \Spec(F) \defeq \Hom_{\text{alg}}(F,\fld) $, is given by
\begin{align*}
  \Spec(F) = \{\lambda_I \mid I \in \I\}
\end{align*}
where, for $ S \in \I $ we have,
\begin{align} \label{eq:fusion-ring-reg-rep}
  \lambda_I(S) =
  \begin{array}{c}
    \begin{tikzpicture}[scale=0.25,every node/.style={inner sep=0,outer sep=-1},yscale=-1]
      \node (v1) at (6,-3.5) {};
      \node (v2) at (6,-6) {};
      \node (v3) at (6,-8.5) {};
      \draw [very thick] (v1) edge (v2);
      \draw [line width=0.2cm,white] (6,-6) ellipse (1.5 and 1.5);
      \draw [thick] (6,-6) ellipse (1.5 and 1.5);
      \draw [line width=0.2cm,white] (v2) edge (v3);
      \draw [very thick] (v2) edge (v3);
      \node at (5.5,-9) {$ I^\vee $};
      \node at (8.5,-6) {$ S $};
    \end{tikzpicture}
  \end{array}
  \in \End_{\C}(I^\vee) = \fld.
\end{align}
We note that taking the trace of both sides of \eqref{eq:fusion-ring-reg-rep} implies that $ \lambda_I(S)d(I) = \Sc_{IS} $. Further still, an analogous argument for $ \lambda_S(I) $ shows that
\begin{align} \label{eq:spectrum-in-terms-of-s-matrix}
  \lambda_I(S)d(I)= \Sc_{IS} =  \lambda_S(I) d(S).
\end{align}
We now consider the inner product on $ \Hom(F,\fld) $ given by
\begin{align*}
  \< a,b \> = \sum_S a({S})\cdot b({S^\vee}).
\end{align*}
By~\cite[Lemma 8.14.1]{Etingof15} $ \Spec(F) $ is orthogonal with respect to this inner product. Therefore for any $ \lambda \in \Spec(F) $ rewriting the corresponding idempotent $ e_\lambda \in F $ in the basis of simple objects gives
\begin{align*}
  e_\lambda = \frac{1}{\|\lambda\|^2} \sum \lambda(S^\vee) \cdot S.
\end{align*}
Indeed, for any $ \mu \in \Spec(F) $ this satisfies $ \mu(e_\lambda)=\delta_{\mu,\lambda} $.

An alternative perspective on $ F $ is provided by the tube category of $ \C $. As described in Remark~\ref{rem:isom_kc_endtid}, $ \End_{\TC}(\tid)=F $. A complete set of orthogonal primitive idempotents in this endomorphism space was identified in~\cite{hardiman_king}, they are given by
\begin{align*}
  \tid_I = \frac{d(I)}{d(\C)} \sum_S d(S)
  \begin{array}{c}
    \begin{tikzpicture}[scale=0.2,every node/.style={inner sep=0,outer sep=-1}]
      \node (v1) at (0,5) {};
      \node (v4) at (0,-5) {};
      \node (v2) at (5,0) {};
      \node (v3) at (-5,0) {};
      \node (v5) at (-2.5,2.5) {};
      \node (v13) at (0,0) {};
      \node (v8) at (2.5,-2.5) {};
      \draw [thick] (v8) edge (v13);
      \draw [line width=0.2cm,white] (v13) ellipse (1.5 and 1.5);
      \draw [thick] (v13) ellipse (1.5 and 1.5);
      \draw [line width=0.2cm,white] (v13) edge (v5);
      \draw [thick] (v13) edge (v5);
      \node at (-3.5,3.5) {$ S $};
      \node at (3.5,-3.5) {$ S $};
      \draw[very thick, red]  (v1) edge (v3);
      \draw[very thick, red]  (v2) edge (v4);
      \draw[very thick]  (v1) edge (v2);
      \draw[very thick]  (v3) edge (v4);
      \node at (-2.5,0) {$ I $};
    \end{tikzpicture}
  \end{array}
\end{align*}
for $ I \in \I $.

The non-degeneracy of the braiding on $\C$ also gives rise to a complete set of primitive idempotents in $\TC$, which provide an alternative labelling of the simple objects in $\RTC$.

\begin{DEF}
  \label{def:primitive-idem-tc}
  For $X$ and $Y$ in $\C$, we consider the following morphism in $\TC$,
  \begin{align*}
    \e_X^Y = \frac{1}{d(\C)} \bigoplus\limits_S d(S)
    \begin{array}{c}
      \begin{tikzpicture}[scale=0.2,every node/.style={inner sep=0,outer sep=-1}]
        \node (v1) at (0,5) {};
        \node (v4) at (0,-5) {};
        \node (v2) at (5,0) {};
        \node (v3) at (-5,0) {};
        \node (v5) at (-2.5,2.5) {};
        \node (v6) at (2.5,-2.5) {};
        \node (v7) at (1.5,3.5) {};
        \node (v11) at (3.5,1.5) {};
        \node (v12) at (-1.5,-3.5) {};
        \node (v10) at (-3.5,-1.5) {};
        \node (v13) at (2.5,-2.5) {};
        \draw [thick] (v7) edge (v10);
        \draw [line width =0.5em,white] (v5) edge (v13);
        \draw [thick] (v5) edge (v13);
        \draw [line width =0.5em,white] (v11) edge (v12);
        \draw [thick] (v11) edge (v12);
        \node at (2.75,4.75) {$ X $};
        \node at (4.75,2.75) {$ Y $};
        \node at (-3.75,3.75) {$ S $};
        \node at (3.75,-3.5) {$ S $};
        \draw[very thick, red]  (v1) edge (v3);
        \draw[very thick, red]  (v2) edge (v4);
        \draw[very thick]  (v1) edge (v2);
        \draw[very thick]  (v3) edge (v4);
      \end{tikzpicture}
    \end{array}
    \in \End_{\TC}(XY).
  \end{align*}
  For any $X$ and $Y$ in $\C$, $\e_{X}^{Y}$ is an idempotent and the set
  \begin{align*}
    \left\{\left(IJ,\e_{I}^{J}\right)\Sharp\right\}_{I,J \in \I}
  \end{align*}
  forms a complete set of simple objects in the Karoubi envelope of $\TC$~\cite{hardiman19a}.
\end{DEF}

\begin{REM}
  \label{rem:equiv-idem}
  To relate the idempotents $\tid_{I}$ to those of Definition~\ref{def:primitive-idem-tc}, we note that $\left(\tid,\tid_{I}\right)\Sharp$ and $\big(II^{\vee},\e_{I}^{I^{\vee}}\big)\Sharp$ are isomorphic as objects in $\RTC$~\cite[Theorem 5.6 and Corollary 5.7]{hardiman_king}.
\end{REM}

The following proposition forms a bridge between these two perspectives on the primitive idempotents in $F$.

\begin{PROP} \label{prop:idempotents-agree}
  Let $ I $ be in $ \I $. Then
  \begin{align} \label{eq:idempotents-agree}
    e_{\lambda_I}=\tid_I.
  \end{align}
  \proof

  Using the same reasoning as was used to derive~\eqref{eq:spectrum-in-terms-of-s-matrix} and then the identification of $ \End_{\TC}(\tid) $ with $ F $, we can rewrite $\tid_I $ as
  \begin{align*}
    \tid_I &= \frac{d(I)}{d(\C)} \sum_S d(S) \lambda_{S^\vee}(I) \ S \\
           &= \frac{d(I)^2}{d(\C)} \sum_S \lambda_{I}(S^\vee) \ S
  \end{align*}
  Therefore,~\eqref{eq:idempotents-agree} is satisfied if and only if  $\|\lambda_I\|^2 = \frac{d(\C)}{d(I)^2}$. We calculate,
  \begin{align*}
    \|\lambda_I\|^2 &= \sum_T \lambda_I(T^\vee) \lambda_I(T) \\
                    &=\frac{1}{d(I)^2} \sum_T \Sc_{IT} \Sc_{IT^\vee}\\
                    &=\frac{1}{d(I)^2} \sum_T \Sc_{IT} \Sc_{TI^\vee} = \frac{(\Sc^2)_{II^\vee}}{d(I)^2}  = \frac{d(\C)}{d(I)^2}
  \end{align*}
  which allows us to conclude. \endproof
\end{PROP}

\section{Introducing $\TM$}
\label{sec:intro-tm}

As before, $\C$ is a modular tensor category. The condition that a module category over $\C$ be pivotal (in the sense of Definition~\ref{def:pivotal-module-category}) has previously appeared as a relevant consideration with respect to the $\TM$ construction~\cite{Hardiman:2019utc}. We now provide a brief recap of the $\TM$ construction.

Let $\D$ be a pivotal monoidal category and let $\M \colon \C \to \D$ be a pivotal functor. We consider the following (contravariant) functor,
\begin{align*}
  \Tr \circ \M \colon \C &\to \Vect \\
  X &\mapsto \Hom_{\D}(\M(X),\tid).
\end{align*}
$\Tr \circ \M$ may then be extended to a functor on $\TC$ by setting, for $\alpha_{G} \in \Hom_{\TC}(X,Y)$,
\begin{align*}
  \TM(\alpha_G) \colon \Hom_{\D}(\M(Y),\tid) &\to \Hom_{\D}(\M(X),\tid)\\
  \beta &\mapsto
          \begin{array}{c}
            \begin{tikzpicture}[scale=0.2,every node/.style={inner sep=0,outer sep=-1},yscale=1]
              \draw [thick, blue] (-4.5,-5) ellipse (2.7 and 3);
              \node [draw,outer sep=0,inner sep=0,minimum size=12, minimum height=13, fill=white] (v1) at (-4.5,-6) {$ \beta $};
              \node [blue,draw,diamond,outer sep=0,inner sep=2,minimum size=18,fill=white] (v2) at (-2.5,-3) {$\alpha $};
              \draw [thick, blue] (v1) to[out=90,in=-135] (v2);
              \node (v3) at (0,0.5) {};
              \draw [thick, blue] (v2) to[out=45,in=-90] (v3);
              \node [blue] at (-0.5,-4.5) {$ G $};
              \node [blue] at (-1.25,0) {$ X $};
            \end{tikzpicture}
          \end{array}
          .
\end{align*}
Decomposing $\TM$ into its simple multiplicities spaces,
\begin{align*}
  \TM_{I}^{J} \defeq \Hom_{\RTC}\left((IJ,\e_{I}^{J})\Sharp, \TM \right)
\end{align*}
allows us to consider the non-negative integer $\I \times \I$-matrix $Z$, with entries
\begin{align*}
  Z_{IJ}\defeq \dim \TM_{I}^{J}.
\end{align*}
As described in the introduction, the equivalence between $ \RTC $ and $ \cZ(\C) $ allows us to induce a spherical monoidal structure on $ \RTC $ which then leads to a corresponding quantum dimension $ d \colon \K(\RTC) \to \fld $. In~\cite{Hardiman:2019utc} it is shown that, under the condition that $\M$ is indecomposable and $d(\TM) = d(\C)$, the matrix $Z$ is a \emph{modular invariant}, i.e.\ a non-negative integer $\I \times \I$-matrix which commutes with the modular data of $\C$ and satisfies $Z_{\tid,\tid} = 1$. In fact, the stronger fact that $\TM$ is a \emph{modular invariant algebra}, in the sense of~\cite[Section 6]{MR2430629}, is proved; this property is revisited in Section~\ref{sec:cft}.

\begin{REM}
  \label{rem:encir-rep-mod-inv-diag}
  As discussed in the introduction, the endomorphism algebra $\End_{\TC}(\tid)$ is canonically isomorphic to the fusion algebra $F$. Therefore, any representation of the tube algebra automatically defines an $F$-module. By comparing the action of $\End_{\TC}(\tid) \cong F$ on $\TM(\tid)$ with the definition of $\Enc$~\eqref{eq:enc_mod}, one sees that $\TM(\tid)$, considered as an $F$-module, is precisely the encircling module $\Enc$ defined in Section~\ref{sec:mod-cat-fusion-alg}.
\end{REM}

Combining this with Theorem~\ref{thm:encircling-module-iso-nim-rep}, we obtain the main result of this paper: the identification of the diagonal part of $\TM$ with the NIM-rep. This result may be considered a categorical generalisation of~\cite[Theorem 4.16, (2)]{Boeckenhauer1999ChiralSO}.

\begin{COR}
  \label{cor:diag_tm_nim_rep}
  Let $\M \colon \C \to \End(\B)$ be a pivotal module category. Then $\TM(\tid)$ is isomorphic to the NIM-rep $\Nim$ as an $F$-module. In particular, for $I \in \I$, the multiplicity of $\lambda_{I} \in \Spec(F)$ within $\Nim$ is equal to $Z_{I,I^{\vee}} = \dim \TM_{I}^{I^{\vee}}$.

  \proof By Remark~\ref{rem:encir-rep-mod-inv-diag}, $\TM(\tid)$ is isomorphic to $\Enc$ as an $F$-module, which, by Theorem~\ref{thm:encircling-module-iso-nim-rep}, is isomorphic to $\Nim$. The second statement follows by noting that, by Proposition~\ref{prop:idempotents-agree}, the $\lambda_{I}$-eigenspace of $\TM(\tid)$ is $\Ima \TM(\tid_{I})$, which, by Remark~\ref{rem:equiv-idem}, is precisely the multiplicity space $\TM_{I}^{I^{\vee}}$. \endproof
\end{COR}

\section{The dimension of $\TM$}

\begin{DEF} \label{def:flat}
  For an object $\F$ in $\RTC$, we define $\F\Flat \in \Obj(\C)$ by
  \begin{align*}
    \F\Flat = \sum_S \F(S) \vprod S,
  \end{align*}
  where the sum ranges over a complete set of simple objects in $\C$. In other words, $\F\Flat$ is the object in $\C$ representing the underlying functor obtained by restricting $\F$ to $\C$-morphisms.
\end{DEF}

The goal of this section is to use Theorem~\ref{thm:encircling-module-iso-nim-rep} to prove that, when $\M$ is an indecomposable pivotal module category $\M \colon \C \to \End(\B)$, the condition $d(\TM) = d(\C)$ is automatically satisfied. The following lemma gives an initial strategy for computing the dimension of objects in $\RTC$.

\begin{LEMMA}
  Let $ \F $ be an object in $ \RTC $. Then
  \begin{align*}
    d(\F) = d(\F\Flat)
  \end{align*}
  where $\F\Flat$ is as in Definition~\ref{def:flat}.

  \proof

  Abbreviating the object $(IJ,\e_{I}^{J})\Sharp$ of $\RTC$ to $\rtcs{I}{J}$, we start by establishing the fact that $ d(\rtcs{I}{J}) = d(IJ) $. By definition of the quantum dimension, we have
  \begin{align*}
    d(\rtcs{I}{J}) \id_{\rtcs{\tid}{\tid}} = \ann(\rtcs{I}{J}) \circ \cre(\rtcs{I}{J})
  \end{align*}
  where $ \cre $ and $ \ann $ denote the creation and annihilation morphisms respectively. As described in the latter part of Section 5 in~\cite{Hardiman:2019utc}, this $ \RTC $-morphism may be described by the $ \TC $-morphism
  \begin{align} \label{eq:tc-quant-dim}
    \begin{array}{c}
      \begin{tikzpicture}[scale=0.2,every node/.style={inner sep=0,outer sep=-1}]
        \node (v1) at (0,5) {};
        \node (v4) at (0,-5) {};
        \node (v2) at (5,0) {};
        \node (v3) at (-5,0) {};
        \node (v5) at (-2.5,2.5) {};
        \node (v6) at (-1.25,1.25) {};
        \node (v13) at (0,0) {};
        \node (v14) at (1.25,-1.25) {};
        \node (v8) at (2.5,-2.5) {};
        \draw [line width=0.15cm,white] (v13) ellipse (2.5 and 2.5);
        \draw [thick] (v13) ellipse (2.5 and 2.5);
        \draw[line width=0.15cm,white]  (v5) edge (v6);
        \draw[thick]  (v5) edge (v6);
        \draw[line width=0.15cm,white]  (v8) edge (v13);
        \draw [thick] (v8) edge (v13);
        \draw [thick] (v14) edge (v6);
        \draw [line width=0.15cm,white] (v13) ellipse (1 and 1);
        \draw [thick] (v13) ellipse (1 and 1);
        \node at (-3.5,3.5) {$ S $};
        \node at (3.5,-3.5) {$ S $};
        \draw[very thick, red]  (v1) edge (v3);
        \draw[very thick, red]  (v2) edge (v4);
        \draw[very thick]  (v1) edge (v2);
        \draw[very thick]  (v3) edge (v4);
        \node at (-1.625,0) {$ \scriptstyle I $};
        \node at (-3.375,0) {$ \scriptstyle J $};
      \end{tikzpicture}
    \end{array}
  \end{align}
  considered as an endomorphism of the idempotent $ \rtcs{\tid}{\tid} $ in the Karoubi envelope. As~\eqref{eq:tc-quant-dim} evaluates to $ d(IJ) \id_{\rtcs{\tid}{\tid}} $, we have proven that $ d(\rtcs{I}{J})=d(IJ) $. To conclude we may now simply compute the dimension of $ \F $ as an object in $ \RTC $,
  \begin{align*}
    d(\F) &= \sum_{I,J} \dim \F_I^J \ d(\rtcs{I}{J}) = \sum_{I,J} \dim \F_I^J \ d(IJ) \\
         &=\sum_{S} \sum_{I,J} \dim \F_I^J \hom_{\C}(S,IJ) \ d(S)\\
         &=\sum_{S} \sum_{I,J} \dim \F_I^J \hom_{\TC}(S,\rtcs{I}{J}) \ d(S)\\
         &= \sum_S \dim \F(S) \ d(S) = d(\F\Flat),
  \end{align*}
  where $\F_{I}^{J} = \Hom_{\RTC}(\rtcs{I}{J},\F)$ and $\hom_{\TC}(S,\rtcs{I}{J})$ denotes the dimension of $\rtcs{I}{J}(S)$.
  \endproof
\end{LEMMA}

\begin{REM}
  Under the equivalence between $\RTC$ and the Drinfeld centre $\cZ(\C)$~\cite{MR1990929, hardiman19a} the fact that $ d(\F)=d(\F\Flat) $ corresponds to the fact that, for any object $ (X,\tau) $ in $ \cZ(\C) $ we have $ d((X,\tau))=d(X) $.
\end{REM}

We note that, $ \F\Flat $ is precisely the object representing the functor one obtains by restricting $ \F $ to $ \C $-morphisms. In particular,
\begin{align*}
  d(\TM) = d(\TM\Flat) = \sum_S \dim \Tr \circ \M (S) \ d(S).
\end{align*}
Furthermore, just as we decategorified $\M$ to obtain the NIM-rep, we can also decategorify  $\Tr \circ \M$:
\begin{align*}
  \begin{array}{rccccc}
    \Tr \circ \M \colon & \C &  \xrightarrow{\M} &  \End(\B) &  \xrightarrow{\Tr} &  \Vect \\
    \chi_{\M}\colon& \K(\C)& \xrightarrow{\phantom{\M}} &  \End(\Nim) & \xrightarrow{\tr} &  \Z.
  \end{array}
\end{align*}
By semisimplicity of $ F $, the ``complexification'' of $ \chi_{\M} $ (which, by abuse of notation we also denote by $ \chi_{\M} $) decomposes into a non-negative integer linear combination of points in $ \Spec(F) $:
\begin{align} \label{eq:describe-multiplicites}
  \chi_{\M} = \sum_{\lambda \in \Spec(F)} [\lambda,\chi_{\M}] \ \lambda,
\end{align}
where $ [\lambda, \chi_{\M}] $ is the dimension of the $ \lambda $-eigenspace of the NIM-rep. As $ \dim \Tr \circ \M(S) = \chi_{\M}({S}) $, we have,
\begin{align*}
  d(\TM) = \sum_S \chi({S})\cdot d(S) = \sum_S \chi({S})\cdot d({S^\vee}) = \<\chi_{\M},d\>.
\end{align*}
Finally, by the orthogonality of $ \Spec F $,
\begin{align*}
  \<\chi_{\M},d\>=\mult\<d,d\>=\mult d(\C).
\end{align*}
Recall that, by definition, $\mult$ is the dimension of the $d$-eigenspace in the NIM-rep $\Nim$. Since $d = \lambda_{\tid}$ (take $I = \tid$ in~\eqref{eq:fusion-ring-reg-rep}), Corollary~\ref{cor:diag_tm_nim_rep} with $I = \tid$ shows that this dimension equals $\dim \TM_{\tid}^{\tid}$. Taking all of this together with Proposition 6.3 in~\cite{Hardiman:2019utc}, we obtain the following.
\begin{PROP}
  \label{prop:dimension-tm}
  Let $\M$ be a pivotal module category over a modular tensor category $\C$. The following are equivalent,
  \begin{enumerate}
  \item $\M$ is indecomposable, \label{cond:indecomp}
  \item $\TM_{\tid}^{\tid} = \fld$, \label{cond:haploid}
  \item $d(\TM) = d(\C)$. \label{cond:dimC}
  \end{enumerate}
  \proof Proposition 6.3 in~\cite{Hardiman:2019utc} proves that \eqref{cond:indecomp} and \eqref{cond:haploid} are equivalent. By the discussion above we have
  \begin{align*}
    d(\TM) = \<\chi_{\M},d\>=\dim \TM_{\tid}^{\tid} \ d(\C),
  \end{align*}
  which, as $d(\C) \neq 0$ (Section 4) and $\fld$ has characteristic zero, proves the equivalence between \eqref{cond:haploid} and \eqref{cond:dimC}. \endproof
\end{PROP}

\begin{COR}
  \label{cor:TM-mod-inv}
  Let $\M$ be an indecomposable pivotal module category over $\C$. Then the matrix $Z$ associated to $\TM$ is a modular invariant, i.e.\ a non-negative integer $\I \times \I$-matrix which commutes with the \emph{modular data} of $\C$ and satisfies $Z_{\tid,\tid} = 1$.
  \proof The result follows directly from~\cite[Corollary 6.9]{Hardiman:2019utc} and Proposition~\ref{prop:dimension-tm}. \endproof
\end{COR}

\section{$\TM$ and the full centre} \label{sec:cft}

In this section we identify $\TM$ with the \emph{full centre} construction~\cite{MR2443266,MR2658689,MR3406516}, thereby connecting the $\TM$ construction to existing literature on algebra objects in the Drinfeld centre (see also~\cite{MR1940282,MR2551797,FuchsSY2023,Yang2022}). As before, $\C$ is a modular tensor category. We suppress the associativity and unit constraints of $\C$ and the monoidal structure of $\M$.

Let $\M \colon \C \to \End(\B)$ be a finite $\C$-module category and let $(X,\sigma)$ be an object of $Z(\C)$, with half-braiding $\sigma_{G} \colon XG \to GX$. Via $\M(\sigma_{G})$, the functor $\M(X) = X \otimes \Blank$ is a $\C$-module functor. A natural transformation $\phi \colon \M(X) \Rightarrow \id_{\B}$ is called a \emph{$\C$-module natural transformation} if
\begin{align}
  \label{eq:module-nat}
  \phi \otimes \id_{\M(G)} = (\id_{\M(G)} \otimes \phi) \circ \M(\sigma_{G}) \colon \M(X)\M(G) \Rightarrow \M(G)
\end{align}
for all $G \in \C$. In components, this says that
\begin{align*}
  \phi_{G \otimes m} = (\id_{G} \otimes \phi_{m}) \circ (\sigma_{G} \otimes \id_{m}) \qquad \text{for all } G \in \C \text{ and } m \in \B.
\end{align*}
We write $\Nat_{\C}(X \otimes \Blank, \id_{\B})$ for the space of $\C$-module natural transformations.

The \emph{full centre} of $\M$ is the commutative algebra $Z(\M) \in Z(\C)$ defined as the internal endomorphism object $\underline{\Hom}(\id_{\B}, \id_{\B})$~\cite[Definition~3.3 and~(3.7)]{MR3406516}. It is characterised by the bijection
\begin{align}
  \label{eq:full-centre-univ-prop}
  \Hom_{Z(\C)}\!\big((X,\sigma),\, Z(\M)\big) \;\cong\; \Nat_{\C}\!\big(X \otimes \Blank,\, \id_{\B}\big),
\end{align}
natural in $(X,\sigma) \in Z(\C)$~\cite[(3.8)]{MR3406516}. Under~\eqref{eq:full-centre-univ-prop}, the unit of $Z(\M)$ corresponds to $\id_{\id_{\B}}$, and its multiplication corresponds to horizontal composition~\cite[Remark~2.9]{MR3406516},
\begin{align*}
  (\phi,\psi) \mapsto \phi \otimes_{\D} \psi = \phi \circ (\id_{\M(X)} \otimes \psi).
\end{align*}

To compare $Z(\M)$ with $\TM$, we need the equivalence $\RTC \simeq Z(\C)$ of Remark~\ref{rem:equiv_zc} in an explicit form, which we now recall.

\begin{REM}
  \label{rem:explicit-equivalence}
  For $G \in \C$, we write
  \begin{align*}
    \widetilde{\cre}_{G} &= (\id_{G^{\vee}} \otimes p_{G}^{-1}) \circ \cre_{G^{\vee}} \colon \tid \to G^{\vee}G, \\
    \widetilde{\ann}_{G} &= \ann_{G^{\vee}} \circ (p_{G} \otimes \id_{G^{\vee}}) \colon GG^{\vee} \to \tid.
  \end{align*}
  For $G,X \in \C$, let $c_{G,X} \in \Hom_{\TC}(XG,GX)$ denote the morphism $\alpha_{G}$ of~\eqref{eq:tc_morphism} with $\alpha = \id_{GXG}$. Every morphism of $\TC$ is a composite of morphisms of $\C$ and morphisms of the form $c_{G,X}$. Indeed, for $\alpha \in \Hom_{\C}(GX,YG)$, the composition rule~\eqref{eq:composition_in_tc} and the zig-zag identity for $(\widetilde{\cre}_{G},\widetilde{\ann}_{G})$ give
  \begin{align*}
    \alpha_{G} = (\widetilde{\ann}_{G} \otimes \id_{Y}) \circ c_{G,G^{\vee}Y} \circ \big( (\id_{G^{\vee}} \otimes \alpha) \circ (\widetilde{\cre}_{G} \otimes \id_{X}) \big).
  \end{align*}
  Hence, by~\cite[Proposition~6.1]{hardiman19a}, an object $\F$ of $\RTC$ is determined by its restriction to $\C$ together with the isomorphisms
  \begin{align*}
    \kappa^{\F}_{G,X} \defeq \F(c_{G,X}) \colon \F(GX) \to \F(XG).
  \end{align*}
  Likewise, a morphism in $\RTC$ is a natural transformation between the restrictions to $\C$ which commutes with these isomorphisms~\cite[Proposition~6.3]{hardiman19a}.

  Fix an isomorphism $\F|_{\C} \cong \Hom_{\C}(\Blank, \F\Flat)$ (Definition~\ref{def:flat}). The equivalence sends $\F$ to $(\F\Flat,\sigma)$, where the half-braiding $\sigma_{G} \colon \F\Flat G \to G \F\Flat$ is determined by
  \begin{align}
    \label{eq:kappa-sigma}
    \kappa^{\F}_{G,X}(f) = (\ann_{G} \otimes \id_{\F\Flat}) \circ (\id_{G^{\vee}} \otimes \sigma_{G}) \circ (\id_{G^{\vee}} \otimes f \otimes \id_{G}) \circ (\widetilde{\cre}_{G} \otimes \id_{XG})
  \end{align}
  for $f \in \Hom_{\C}(GX, \F\Flat)$. In words, $\kappa^{\F}_{G,X}$ carries the strand $G$ once around the cylinder, using $\sigma_{G}$ to cross $\F\Flat$. In particular, $\Blank\Flat$ corresponds to the forgetful functor $Z(\C) \to \C$.

  In~\cite{hardiman19a}, half-braidings are written $\tau_{G} \colon G\F\Flat \to \F\Flat G$. Here we use $\sigma_{G} = \tau_{G}^{-1}$, which matches the conventions of~\cite{MR3406516}; the map $(X,\tau) \mapsto (X,\tau^{-1})$ is an isomorphism of monoidal categories between the two versions of $Z(\C)$. In terms of $\tau$, the right-hand side of~\eqref{eq:kappa-sigma} equals
  \begin{align*}
    (\id_{\F\Flat} \otimes \ann_{G}) \circ (\tau_{G^{\vee}} \otimes \id_{G}) \circ (\id_{G^{\vee}} \otimes f \otimes \id_{G}) \circ (\widetilde{\cre}_{G} \otimes \id_{XG}).
  \end{align*}
  This follows from the half-braiding axiom $\tau_{G^{\vee}G} = (\tau_{G^{\vee}} \otimes \id_{G}) \circ (\id_{G^{\vee}} \otimes \tau_{G})$ together with naturality of $\tau$ with respect to $\ann_{G}$. It is the formula of~\cite[Section~7]{hardiman19a} for representable functors, with $\alpha = \id_{GXG}$.

  Finally, the equivalence is monoidal, where $(X,e)\Sharp \otimes (Y,e')\Sharp = (XY, e \otimes_{\TC} e')\Sharp$~\cite[Sections~5 and~6]{Hardiman:2019utc}. For $X \in \C$, let $u_{X} \colon X \to (X\Sharp)\Flat$ be the morphism corresponding to $\id_{X} \in X\Sharp(X)$. Then
  \begin{align*}
    u_{XY} = u_{X} \otimes u_{Y} \qquad \text{under } (XY)\Sharp \cong X\Sharp \otimes Y\Sharp.
  \end{align*}
\end{REM}

\begin{THM}
  \label{thm:tm-iso-full-centre}
  Let $\M \colon \C \to \End(\B)$ be a finite pivotal module category over a modular tensor category $\C$. For every $V \in \RTC$, corresponding to $(V\Flat,\sigma) \in Z(\C)$, restriction to $\C$ induces a bijection
  \begin{align}
    \label{eq:tm-univ-prop}
    \Psi_{V} \colon \Hom_{\RTC}(V, \TM) \xrightarrow{\ \sim\ } \Nat_{\C}(V\Flat \otimes \Blank, \id_{\B}),
  \end{align}
  natural in $V$. Consequently, under the equivalence $\RTC \simeq Z(\C)$, there is a canonical isomorphism of algebras
  \begin{align*}
    \TM \cong Z(\M).
  \end{align*}
  \proof
  Let $V \in \RTC$. Restriction to $\C$ is injective on $\Hom_{\RTC}(V,\TM)$, since $\C$ and $\TC$ have the same objects. By the Yoneda lemma, a natural transformation $\eta \colon V|_{\C} \Rightarrow \TM|_{\C}$ corresponds to the element
  \begin{align*}
    \phi = \eta_{V\Flat}(\id_{V\Flat}) \in \TM(V\Flat) = \Nat(\M(V\Flat), \id_{\B}),
  \end{align*}
  with $\eta_{X}(f) = \phi \circ \M(f)$ for $f \in \Hom_{\C}(X, V\Flat)$. We show that $\eta$ commutes with the isomorphisms $\kappa$ if and only if $\phi$ satisfies~\eqref{eq:module-nat}. This proves~\eqref{eq:tm-univ-prop}, setting $\Psi_{V}(\eta) = \phi$. Naturality in $V$ is clear, as both sides are functorial by precomposition.

  Specialising the definition of $\TM$ in Section~\ref{sec:intro-tm} to $\alpha = \id_{GXG}$, and using that the pivotal structure on $\D$ restricts to $\M(p)$ on the image of $\M$, we have
  \begin{align*}
    \kappa^{\TM}_{G,X}(\beta) = \M(\ann_{G}) \circ \big( \id_{\M(G^{\vee})} \otimes \beta \otimes \id_{\M(G)} \big) \circ \big( \M(\widetilde{\cre}_{G}) \otimes \id_{\M(XG)} \big)
  \end{align*}
  for $\beta \in \TM(GX)$. Thus $\eta$ commutes with $\kappa$ if and only if
  \begin{align}
    \label{eq:kappa-compat}
    \begin{split}
      \phi \circ \M\big(\kappa^{V}_{G,X}(f)\big) = {} & \M(\ann_{G}) \circ \big( \id_{\M(G^{\vee})} \otimes (\phi \circ \M(f)) \otimes \id_{\M(G)} \big) \\
      & \circ \big( \M(\widetilde{\cre}_{G}) \otimes \id_{\M(XG)} \big)
    \end{split}
  \end{align}
  for all $G, X \in \C$ and $f \in \Hom_{\C}(GX, V\Flat)$.

  Both sides of~\eqref{eq:kappa-compat} are natural in $X$. Moreover, every $f$ factors as
  \begin{align*}
    f = f_{0} \circ (\id_{G} \otimes \hat{f}), \qquad f_{0} = \widetilde{\ann}_{G} \otimes \id_{V\Flat} \colon GG^{\vee}V\Flat \to V\Flat, \quad \hat{f} \colon X \to G^{\vee}V\Flat.
  \end{align*}
  So it suffices to consider $X = G^{\vee}V\Flat$ and $f = f_{0}$. For this choice, the zig-zag identity $(\id_{G^{\vee}} \otimes \widetilde{\ann}_{G}) \circ (\widetilde{\cre}_{G} \otimes \id_{G^{\vee}}) = \id_{G^{\vee}}$ simplifies both~\eqref{eq:kappa-sigma} and the right-hand side of~\eqref{eq:kappa-compat}. Equation~\eqref{eq:kappa-compat} then becomes
  \begin{align*}
    & \M(\ann_{G}) \circ \Big( \id_{\M(G^{\vee})} \otimes \big( (\id_{\M(G)} \otimes \phi) \circ \M(\sigma_{G}) \big) \Big) \\
    & \qquad = \M(\ann_{G}) \circ \Big( \id_{\M(G^{\vee})} \otimes \big( \phi \otimes \id_{\M(G)} \big) \Big)
  \end{align*}
  as morphisms $\M(G^{\vee})\M(V\Flat)\M(G) \to \id_{\B}$. Since $\M(G^{\vee})$ is a left dual of $\M(G)$ in $\D$, the map
  \begin{align*}
    \Hom_{\D}(Q, \M(G)) \to \Hom_{\D}(\M(G^{\vee})Q, \tid), \qquad \psi \mapsto \M(\ann_{G}) \circ (\id_{\M(G^{\vee})} \otimes \psi),
  \end{align*}
  is a bijection for every $Q \in \D$. So this equation is equivalent to~\eqref{eq:module-nat}.

  By~\eqref{eq:tm-univ-prop},~\eqref{eq:full-centre-univ-prop} and the Yoneda lemma, there is a unique isomorphism $\TM \cong Z(\M)$ compatible with these bijections. It remains to compare the algebra structures. The unit of $\TM$ is~\cite[Proposition~6.2]{Hardiman:2019utc}
  \begin{align*}
    u = \id_{\tid_{\D}} \in \TM_{\tid}^{\tid} = \Hom_{\RTC}\big((\tid,\e_{\tid}^{\tid})\Sharp, \TM\big),
  \end{align*}
  and $\Psi_{(\tid,\e_{\tid}^{\tid})\Sharp}(u) = \id_{\id_{\B}}$ corresponds to the unit of $Z(\M)$. For the multiplication $\nabla$ of $\TM$~\cite[Proposition~6.2]{Hardiman:2019utc}, it suffices to show that
  \begin{align}
    \label{eq:mult-compat}
    \Psi_{V_{1} \otimes V_{2}}\big( \nabla \circ (\eta_{1} \otimes \eta_{2}) \big) = \Psi_{V_{1}}(\eta_{1}) \otimes_{\D} \Psi_{V_{2}}(\eta_{2})
  \end{align}
  for all $\eta_{i} \colon V_{i} \to \TM$. Every object of $\RTC$ is a quotient of a representable functor, so by naturality we may assume $V_{i} = X_{i}\Sharp$.

  Write $\phi_{i} = \Psi_{V_{i}}(\eta_{i})$. By the Yoneda lemma, $\eta_{i}$ corresponds to
  \begin{align*}
    f_{i} = (\eta_{i})_{X_{i}}(\id_{X_{i}}) = \phi_{i} \circ \M(u_{X_{i}}) \in \TM(X_{i}).
  \end{align*}
  By the definition of $\nabla$~\cite[Section~6]{Hardiman:2019utc}, $\nabla \circ (\eta_{1} \otimes \eta_{2}) \colon (X_{1}X_{2})\Sharp \to \TM$ corresponds to $f_{1} \otimes_{\D} f_{2} \in \TM(X_{1}X_{2})$. Using $u_{X_{1}X_{2}} = u_{X_{1}} \otimes u_{X_{2}}$, it follows that both sides of~\eqref{eq:mult-compat}, composed with $\M(u_{X_{1}X_{2}})$, equal $f_{1} \otimes_{\D} f_{2}$.

  Both sides of~\eqref{eq:mult-compat} are $\C$-module natural transformations; for the right-hand side, this is because horizontal composites of $\C$-module natural transformations are $\C$-module natural transformations for the tensor product half-braiding. Composition with $\M(u_{X_{1}X_{2}})$ is injective on $\C$-module natural transformations, because composed with $\Psi^{-1}$ it is the Yoneda bijection. Therefore~\eqref{eq:mult-compat} holds. Since the multiplication of $Z(\M)$ corresponds to $\otimes_{\D}$ under~\eqref{eq:full-centre-univ-prop}, the isomorphism $\TM \cong Z(\M)$ is an isomorphism of algebras.
  \endproof
\end{THM}

\begin{REM}
  \label{rem:regular-module}
  Consider the regular module category $\M \colon \C \to \End(\C)$, $\M(X) = X \otimes \Blank$. A $\C$-module natural transformation $\phi \colon V\Flat \otimes \Blank \Rightarrow \id_{\C}$ is determined by $\phi_{\tid} \in \Hom_{\C}(V\Flat, \tid)$, via
  \begin{align*}
    \phi_{m} = (\id_{m} \otimes \phi_{\tid}) \circ \sigma_{m},
  \end{align*}
  and every element of $\Hom_{\C}(V\Flat,\tid)$ arises in this way. For $V = (IJ,\e_{I}^{J})\Sharp$ we have $V\Flat \cong IJ$ (see the proof of Lemma 6.2), so Theorem~\ref{thm:tm-iso-full-centre} gives
  \begin{align*}
    \dim \TM_{I}^{J} = \dim \Hom_{\C}(IJ, \tid) = \delta_{J,I^{\vee}}.
  \end{align*}
  Thus $\TM \cong \tid\Sharp$ by~\eqref{eq:tid_decomposition}, and $Z$ is the charge-conjugation matrix $Z_{IJ} = \delta_{J,I^{\vee}}$, in agreement with~\cite{MR1940282}.

  This example also illustrates the role of $\C$-module natural transformations. The space $\TM(\tid) = \Hom_{\RTC}(\tid\Sharp,\TM) = \End(\id_{\C})$ has dimension $|\I|$, whereas $\TM_{\tid}^{\tid}$, which corresponds to $\Nat_{\C}(\id_{\C},\id_{\C})$, is one-dimensional.
\end{REM}

We now relate $\TM$ to one of the main applications of the full centre. A rational conformal field theory on surfaces with boundary involves both bulk fields and boundary fields. Kong and Runkel introduced \emph{Cardy $\C$-algebras}~\cite{MR2551797} to capture this algebraically. A Cardy $\C$-algebra consists of an algebra of bulk fields in $Z(\C)$, an algebra of boundary fields in $\C$, and a map between them, subject to compatibility conditions; see Definition~\ref{def:cardy-algebra}. These conditions encode the sewing constraints of the theory in genus zero and one~\cite{MR2551797,KONG2014604}. Kong and Runkel showed that any special symmetric Frobenius algebra $A$ in $\C$, together with its full centre $Z(A)$, gives a Cardy $\C$-algebra~\cite[Theorem~3.18]{MR2551797}.

In the rest of this section we explain how this construction relates to ours. We start from a module category $\M$ rather than an algebra, and use Theorem~\ref{thm:tm-iso-full-centre} to show that $\TM$ is the algebra of bulk fields of such a Cardy $\C$-algebra (Corollary~\ref{cor:tm-cardy-alg}).

For the remainder of this section we take $\fld = \mathbb{C}$. This is needed in order to apply the results of~\cite{MR2551797} and~\cite{MR3019263}: the modular tensor categories of~\cite{MR2551797} are $\mathbb{C}$-linear by definition, and the arguments of~\cite{MR3019263} use complex conjugation (for instance, the dimension matrix of a module category is Hermitian~\cite[Proposition~5.6]{MR3019263}). To use the results of~\cite{MR2551797}, we also need to fix the braided structure on $\RTC$.

\begin{REM}
  \label{rem:braiding-convention}
  In the conventions of Remark~\ref{rem:explicit-equivalence}, with half-braidings $\sigma_{G} \colon XG \to GX$, we equip $Z(\C)$, and hence $\RTC$, with the braiding
  \begin{align*}
    c_{(X,\sigma),(Y,\sigma')} = \sigma_{Y} \colon XY \to YX,
  \end{align*}
  as in~\cite{MR3406516}. Let $\CC = \C^{+} \boxtimes \C^{-}$ carry the braiding $c \boxtimes c^{-1}$, as in~\cite[(2.38)]{MR2551797}. Then the functor
  \begin{align*}
    \CC \to Z(\C), \qquad U \boxtimes V \mapsto \big( UV, (c_{U,\Blank} \otimes \id_{V}) \circ (\id_{U} \otimes c^{-1}_{\Blank,V}) \big)
  \end{align*}
  is a braided equivalence~\cite{Mueger2003}, with monoidal structure given by $\id_{U} \otimes c^{-1}_{U',V} \otimes \id_{V'}$. Composed with the forgetful functor, it is the tensor functor $T \colon \CC \to \C$ of~\cite[Section~2.4 and~(2.40)]{MR2551797}, so $\Blank\Flat$ corresponds to $T$.

  In terms of $\tau = \sigma^{-1}$, this braiding is the reverse of the braiding $\gamma_{X}$ used in~\cite[Section~7]{hardiman19a}. Reversing the braiding does not affect whether an algebra is commutative.
\end{REM}

\begin{DEF}[{\cite[Definition~3.14]{MR2551797}}]
  \label{def:cardy-algebra}
  Let $\C$ be a modular tensor category, and identify $\RTC \simeq Z(\C) \simeq \CC$ as in Remark~\ref{rem:braiding-convention}. A \emph{Cardy $\C$-algebra} is a triple $(\opalg, \clalg, \iota)$ consisting of the following data.
  \begin{itemize}
  \item $\opalg$ is a symmetric Frobenius algebra in $\C$.
  \item $\clalg$ is a commutative symmetric Frobenius algebra in $\RTC$ which is modular invariant in the sense of~\cite[Definition~3.1]{MR2551797} (cf.~\cite[Section~6]{MR2430629}).
  \item $\iota \colon \clalg\Flat \to \opalg$ is an algebra homomorphism satisfying the \emph{centre condition}~\cite[(3.20)]{MR2551797} and the \emph{Cardy condition}~\cite[(3.22)]{MR2551797}.
  \end{itemize}
\end{DEF}

Let $A$ be a haploid special symmetric Frobenius algebra in $\C$ (see~\cite[Section~2.2]{MR2551797}), let $\B = \Mod_{\C}(A)$ be the category of right $A$-modules in $\C$, and let
\begin{align*}
  \M \colon \C \to \End(\B), \qquad X \mapsto X \otimes \Blank,
\end{align*}
be the corresponding module category. As $A$ is special, it is separable, so $\B$ is finite semisimple. By~\cite[Proposition~6.8]{MR3019263}, $\B$ admits a module trace, so $\M$ is pivotal by Proposition~\ref{prop:module-trace} and $\TM$ is defined. Conversely, every indecomposable pivotal module category arises (non-uniquely) in this way; see Corollary~\ref{cor:tm-cardy-alg}.

Recall from~\cite[Section~2.4 and Theorem~2.21]{MR2551797} the functor $R \colon \C \to \CC \simeq \RTC$, which is both left and right adjoint to $T = \Blank\Flat$. For $g \colon V \to R(A)$ we write $g\Natural \colon V\Flat \to A$ for its adjunct, so that
\begin{align*}
  \Hom_{\RTC}(V, R(A)) \to \Hom_{\C}(V\Flat, A), \qquad g \mapsto g\Natural,
\end{align*}
is a bijection natural in $V$. As $R$ is lax monoidal, $R(A)$ is an algebra, with
\begin{align*}
  \big(m_{R(A)} \circ (g_{1} \otimes g_{2})\big)\Natural = m_{A} \circ (g_{1}\Natural \otimes g_{2}\Natural), \qquad \eta_{R(A)}\Natural = \eta_{A}.
\end{align*}

The \emph{full centre} $Z(A)$ of $A$ is the left centre of the algebra $R(A)$~\cite[Definition~3.17]{MR2551797} (see also~\cite[Definition~4.9]{MR2443266}). We write $e_{A} \colon Z(A) \hookrightarrow R(A)$ for the embedding, which is an algebra map~\cite[Lemma~3.19]{MR2551797}. By~\cite[(3.30)]{MR2551797}, $Z(A)$ is the largest subobject $C \subseteq R(A)$ such that
\begin{align*}
  m_{R(A)} \circ c_{R(A),R(A)} = m_{R(A)} \colon C \otimes R(A) \to R(A),
\end{align*}
where $c$ denotes the braiding.

Finally, consider the map $i_{A} \colon \TM\Flat \to A$ whose $S$-component is
\begin{align*}
  \TM(S) = \Nat(S \otimes \Blank, \id_{\B}) &\to \Hom_{\C}(S,A) \\
  \phi &\mapsto \phi_{A} \circ (\id_{S} \otimes \eta_{A}),
\end{align*}
where $\eta_{A} \colon \tid \to A$ is the unit of $A$. Note that $\TM(S)$ consists of \emph{all} natural transformations $S \otimes \Blank \Rightarrow \id_{\B}$, not only the $\C$-module ones.

\begin{PROP}
  \label{prop:full-centre-algebra}
  Let $A$ and $\M$ be as above. Then there is a unique morphism $\hat{e}_{A} \colon \TM \to R(A)$ whose adjunct is $i_{A}$, i.e.\ $\hat{e}_{A}\Natural = i_{A}$. The morphism $\hat{e}_{A}$ is an algebra monomorphism with image $Z(A)$. Hence $\hat{e}_{A} = e_{A} \circ \theta_{A}$, where $\theta_{A} \colon \TM \to Z(A)$ is an isomorphism of algebras.
  \proof
  For $V \in \RTC$, define
  \begin{align*}
    \rho_{V} \colon \Nat_{\C}(V\Flat \otimes \Blank, \id_{\B}) \to \Hom_{\C}(V\Flat, A), \qquad \rho_{V}(\phi) = \phi_{A} \circ (\id_{V\Flat} \otimes \eta_{A}).
  \end{align*}
  This is natural in $V$. So, by Theorem~\ref{thm:tm-iso-full-centre} and the Yoneda lemma, there is a unique morphism $\hat{e}_{A} \colon \TM \to R(A)$ such that
  \begin{align*}
    (\hat{e}_{A} \circ \eta)\Natural = \rho_{V}(\Psi_{V}(\eta)) \qquad \text{for all } \eta \colon V \to \TM.
  \end{align*}
  Let $\varepsilon = \Psi_{\TM}(\id_{\TM})$, so that $\hat{e}_{A}\Natural = \varepsilon_{A} \circ (\id_{\TM\Flat} \otimes \eta_{A})$. For $S \in \C$, an element $\phi \in \TM(S)$ corresponds to $\bar{\phi} \in \Hom_{\C}(S, \TM\Flat)$ with $\phi = \varepsilon \circ \M(\bar{\phi})$. Hence $\hat{e}_{A}\Natural \circ \bar{\phi} = \phi_{A} \circ (\id_{S} \otimes \eta_{A})$, i.e.\ $\hat{e}_{A}\Natural = i_{A}$. Uniqueness is clear, as $g \mapsto g\Natural$ is bijective.

  \emph{$\hat{e}_{A}$ is a monomorphism.} Let $\phi \in \Nat_{\C}(V\Flat \otimes \Blank, \id_{\B})$. As $\phi_{A}$ is a morphism of right $A$-modules, and by~\eqref{eq:module-nat},
  \begin{align*}
    \phi_{A} = m_{A} \circ (\rho_{V}(\phi) \otimes \id_{A}), \qquad \phi_{X \otimes A} = (\id_{X} \otimes \phi_{A}) \circ (\sigma_{X} \otimes \id_{A})
  \end{align*}
  for all $X \in \C$. Every object of $\B$ is a direct summand of a free module $X \otimes A$~\cite[Section~3.1]{MR1976459}, so by naturality $\phi$ is determined by $\rho_{V}(\phi)$. Thus each $\rho_{V}$ is injective, i.e.\ composition with $\hat{e}_{A}$ is injective on every $\Hom$-space.

  \emph{$\hat{e}_{A}$ is an algebra map.} For $\phi_{i} \in \Nat_{\C}(V_{i}\Flat \otimes \Blank, \id_{\B})$, the formula for $(\phi_{1})_{A}$ above gives
  \begin{align*}
    \rho(\phi_{1} \otimes_{\D} \phi_{2}) = (\phi_{1})_{A} \circ \big(\id_{V_{1}\Flat} \otimes (\phi_{2})_{A}\big) \circ (\id \otimes \eta_{A}) = m_{A} \circ \big(\rho(\phi_{1}) \otimes \rho(\phi_{2})\big),
  \end{align*}
  and $\rho(\id_{\id_{\B}}) = \eta_{A}$. By Theorem~\ref{thm:tm-iso-full-centre}, the multiplication and unit of $\TM$ correspond under $\Psi$ to $\otimes_{\D}$ and $\id_{\id_{\B}}$. Those of $R(A)$ correspond under $\Blank\Natural$ to $m_{A} \circ (\Blank \otimes \Blank)$ and $\eta_{A}$. So $\hat{e}_{A}$ is an algebra map.

  \emph{The image of $\hat{e}_{A}$ is $Z(A)$.} Let $\eta \colon V \to \TM$, let $\phi = \Psi_{V}(\eta)$, and let $g \colon W \to R(A)$. By Remark~\ref{rem:braiding-convention}, $c_{V,W} = \sigma_{W\Flat}$, where $\sigma$ is the half-braiding of $V$. Therefore
  \begin{align*}
    \big( m_{R(A)} \circ c_{R(A),R(A)} \circ (\hat{e}_{A}\eta \otimes g) \big)\Natural &= m_{A} \circ \big(g\Natural \otimes \rho(\phi)\big) \circ \sigma_{W\Flat}, \\
    \big( m_{R(A)} \circ (\hat{e}_{A}\eta \otimes g) \big)\Natural &= m_{A} \circ \big(\rho(\phi) \otimes g\Natural\big).
  \end{align*}
  These agree. To see this, apply naturality of $\phi$ to the $A$-module map
  \begin{align*}
    m_{A} \circ (g\Natural \otimes \id_{A}) \colon W\Flat \otimes A \to A,
  \end{align*}
  use the formulae for $\phi_{W\Flat \otimes A}$ and $\phi_{A}$ above, and precompose with $\id \otimes \eta_{A}$. Hence the image of $\hat{e}_{A}$ satisfies the defining property of the left centre, and so it is contained in $Z(A)$.

  On the other hand, $\TM \cong Z(\M)$ by Theorem~\ref{thm:tm-iso-full-centre}, and $Z(\M) \cong Z(A)$~\cite[Section~3.1]{MR3406516} (see also~\cite[Section~8]{MR2658689}). So the image of $\hat{e}_{A}$ is a subobject of $Z(A)$ isomorphic to $Z(A)$. As $\RTC$ is semisimple with finite-dimensional $\Hom$-spaces, it is all of $Z(A)$. Finally, $e_{A}$ is a monomorphism and an algebra map, so $\theta_{A} = e_{A}^{-1} \circ \hat{e}_{A}$ is an isomorphism of algebras.
  \endproof
\end{PROP}

\begin{COR}
  \label{cor:tm-cardy-alg}
  Let $\fld = \mathbb{C}$, let $\C$ be a modular tensor category, and let $\M \colon \C \to \End(\B)$ be an indecomposable pivotal module category over $\C$. Then there exist haploid special symmetric Frobenius algebras $A$ in $\C$ with $\M \simeq \Mod_{\C}(A)$; for instance, $A = \underline{\End}(m)$ for any simple $m \in \B$. For any such $A$, identifying $\M$ with $\Mod_{\C}(A)$, the triple $(A, \TM, i_{A})$ is a Cardy $\C$-algebra, where $\TM$ carries the commutative symmetric Frobenius structure transported from $Z(A)$ along $\theta_{A}$.
  \proof
  Let $m \in \B$ be simple and let $A = \underline{\End}(m)$. By Proposition~\ref{prop:module-trace}, $\B$ admits a module trace, so $A$ is a haploid special symmetric Frobenius algebra by~\cite[Theorem~6.6]{MR3019263}. By~\cite[Theorem~3.1]{MR1976459} and its proof,
  \begin{align*}
    \B \to \Mod_{\C}(A), \qquad n \mapsto \underline{\Hom}(m,n),
  \end{align*}
  is an equivalence of module categories.

  Now let $A$ be any haploid special symmetric Frobenius algebra with $\M \simeq \Mod_{\C}(A)$. By~\cite[Theorem~3.18]{MR2551797}, $(A|Z(A), e_{A})$ is a Cardy $\C|\CC$-algebra in the sense of~\cite[Definition~3.7]{MR2551797}. By~\cite[Lemmas~3.11 and~3.13]{MR2551797}, this means that $(A, Z(A), e_{A}\Natural)$ is a Cardy $\C$-algebra in the sense of Definition~\ref{def:cardy-algebra}. Transporting this structure along the isomorphism of algebras $\theta_{A}$ of Proposition~\ref{prop:full-centre-algebra}, and using
  \begin{align*}
    (e_{A} \circ \theta_{A})\Natural = \hat{e}_{A}\Natural = i_{A},
  \end{align*}
  shows that $(A, \TM, i_{A})$ is a Cardy $\C$-algebra.
  \endproof
\end{COR}

\begin{REM}
  \label{rem:choice-of-algebra}
  The hypotheses on $\M$ in Corollary~\ref{cor:tm-cardy-alg} are exactly what is needed. If $A$ is any haploid special symmetric Frobenius algebra, then $\Mod_{\C}(A)$ is pivotal, as noted above. It is also indecomposable, since $A$ is a simple $A$-module and every $A$-module is a direct summand of a free module $X \otimes A$.

  The algebra $A$ is not unique. Different simple objects $m$ give Morita equivalent, but in general non-isomorphic, algebras. So Corollary~\ref{cor:tm-cardy-alg} yields Cardy $\C$-algebras with the same bulk algebra $\TM$ and different boundary algebras. The Frobenius structure transported to $\TM$ may also depend on $A$, but only up to rescaling its counit and coproduct, since $\TM$ is haploid by Proposition~\ref{prop:dimension-tm}.
\end{REM}

\appendix

\newpage\section{Comparing pivotal structures on a multifusion category} \label{app:compare-pivotal-struct}

Let $\A$ be a multifusion category and let $p_{A},q_{A}\colon A \to A^{\vee\vee}$ be two pivotal structures on $\A$. For $S$ a simple object in $\A$, we consider the scalar
\begin{align*}
  \mu_{S} = p_{S}^{-1} \circ q_{S} \in \End_{\A}(S) = \fld.
\end{align*}
\begin{LEMMA}
  \label{lem:comp-pivot-struct}
  The $\mu_{S}$ numbers satisfy the following properties
  \begin{enumerate}[label=(\roman*)]
  \item $\mu_{S}\mu_{S^{\vee}} = 1$ \label{eq:mu-inverse}
  \item $\mu_{S} = 1$, for all $S < \tid$ \label{eq:mu-tid}
  \item $\mu_{ST} = \mu_{S}\mu_{T}$, for all $S,T$ such that $ST$ is simple\label{eq:mu-tensor}.
  \end{enumerate}
  \proof
  Recall that any pivotal structure $p$ satisfies $p_{X}^{\vee} = p_{X^{\vee}}^{-1}$. Therefore,
  \begin{align*}
    \mu_{S}\id_{S^{\vee}} =
    \begin{array}{c}
      \begin{tikzpicture}[scale=0.5,every node/.style={inner sep=0,outer sep=-1}]
        \node (v1) at (-1,2.75) {};
        \node (v2) at (-1,-0.25) {};
        \node (v3) at (0,-0.25) {};
        \node (v4) at (0,2.25) {};
        \node (v5) at (1,2.25) {};
        \node (v6) at (1,-0.75) {};
        \draw [thick] (v1) edge (v2);
        \draw [thick] (v2) to[out=-90,in=-90] (v3);
        \draw [thick] (v3) edge (v4);
        \draw [thick] (v4) to[out=90,in=90] (v5);
        \draw [thick] (v5) edge (v6);
        \node [draw,outer sep=0,inner sep=0.3,minimum width=17, minimum height = 13, fill=white] at (0,1.625) {$\scriptstyle q_{S}$};
        \node [draw,outer sep=0,inner sep=0.3,minimum width=17, minimum height = 13, fill=white] at (0,0.375) {$\scriptstyle p^{-1}_{S}$};
      \end{tikzpicture}
    \end{array}
    =
    \begin{array}{c}
      \begin{tikzpicture}[scale=0.5,every node/.style={inner sep=0,outer sep=-1}]
        \node (v1) at (-1,2) {};
        \node (v2) at (-1,-0.25) {};
        \node (v3) at (0,-0.25) {};
        \node (v4) at (0,1) {};
        \node (v5) at (1,1) {};
        \node (v6) at (1,-0.25) {};
        \node (v7) at (2,-0.25) {};
        \node (v8) at (2,1) {};
        \node (v9) at (3,1) {};
        \node (v10) at (3,-1.5) {};
        \draw [thick] (v1) edge (v2);
        \draw [thick] (v2) to[out=-90,in=-90] (v3);
        \draw [thick] (v3) edge (v4);
        \draw [thick] (v4) to[out=90,in=90] (v5);
        \draw [thick] (v5) edge (v6);
        \draw [thick] (v6) to[out=-90,in=-90] (v7);
        \draw [thick] (v8) to[out=90,in=90] (v9);
        \draw [thick] (v7) edge (v8);
        \draw [thick] (v9) edge (v10);
        \node [draw,outer sep=0,inner sep=0.3,minimum width=17, minimum height = 13, fill=white] at (2,0.25) {$\scriptstyle q_{S}$};
        \node [draw,outer sep=0,inner sep=0.3,minimum width=17, minimum height = 13, fill=white] at (0,0.375) {$\scriptstyle p^{-1}_{S}$};
      \end{tikzpicture}
    \end{array}
    =
    \begin{array}{c}
      \begin{tikzpicture}[scale=0.5,every node/.style={inner sep=0,outer sep=-1}]
        \node (v1) at (-1,2.75) {};
        \node (v2) at (-1,-0.75) {};
        \draw [thick] (v1) edge (v2);
        \node [draw,outer sep=0,inner sep=0.3,minimum width=18, minimum height = 15, fill=white] at (-1,0.25) {$\scriptstyle q^{-1}_{S^{\vee}}$};
        \node [draw,outer sep=0,inner sep=0.3,minimum width=18, minimum height = 15, fill=white] at (-1,1.75) {$\scriptstyle p_{S^{\vee}}$};
      \end{tikzpicture}
    \end{array}
    = \frac{1}{\mu_{S^{\vee}}} \id_{S^{\vee}},
  \end{align*}
  implying \ref{eq:mu-inverse}. To prove \ref{eq:mu-tid}, we recall that, by \cite[Corollary 4.3.2]{Etingof15},
  \begin{align*}
    \bigoplus_{\substack{S \in \I, \, S < \tid}} S = \tid.
  \end{align*}
  Therefore, for any pivotal structure $p$, $\bigoplus p_{S} = p_{\tid} = \id_{\tid} = \bigoplus \id_{S}$, where the sum ranges over $\{S \in \I \mid S < \tid\}$. Finally, to prove \ref{eq:mu-tensor}, we simply compute,
  \begin{align*}
    \mu_{ST} = p_{ST}^{-1} \circ q_{ST} &= (p_{S}^{-1} \otimes p_{T}^{-1}) \circ (q_{S} \otimes q_{T}) \\
                                        &= (p_{S}^{-1} \circ q_{S}) \otimes (p_{T}^{-1} \circ q_{T}) \\
                                        &= \mu_{S}\mu_{T}.
  \end{align*} \endproof
\end{LEMMA}

\bibliographystyle{alpha}
\bibliography{bibli}

\end{document}